\documentclass[12pt]{article}
\usepackage[american]{babel}
\usepackage{amsmath}
\usepackage{latexsym}
\usepackage{amssymb}
\usepackage{theorem}
\usepackage{diagrams}
\usepackage{graphics}
\usepackage{times}
\usepackage{pstcol}
\usepackage[matrix,poly,arc]{xy}
\theoremstyle{change}
\newtheorem{Thm}{Theorem}[section]
\newtheorem{Cor}[Thm]{Corollary}
\newtheorem{Prop}[Thm]{Proposition}
\newtheorem{Lem}[Thm]{Lemma}
{\theorembodyfont{\rmfamily}
\newtheorem{Num}[Thm]{}

\newtheorem{Def}[Thm]{Definition}}
\newcommand{\proof}{\par\medskip\rm\emph{Proof. }}
\newcommand{\skop}{\par\medskip\rm\emph{Sketch of proof. }}
\newcommand{\qed}{\ \hglue 0pt plus 1filll $\Box$}

\newcommand{\too}{\rTo}
\newcommand{\mapstoo}{\longmapsto}

\renewcommand{\SS}{\mathbb{S}}
\newcommand{\RR}{\mathbb{R}}
\newcommand{\ZZ}{\mathbb{Z}}

\newcommand{\CC}{\mathbb{C}}
\newcommand{\LL}{\mathbb{L}}
\renewcommand{\AA}{\mathbb{A}}

\newcommand{\SKIP}[1]{}
\newcommand{\SO}{\mathrm{SO}}
\newcommand{\SU}{\mathrm{SU}}
\newcommand{\SL}{\mathrm{SL}}
\newcommand{\U}{\mathrm{U}}
\newcommand{\GL}{\mathrm{GL}}
\newcommand{\PGL}{\mathrm{PGL}}

\newcommand{\fsl}{\mathfrak{sl}}

\newcommand{\fsu}{\mathfrak{su}}

\newcommand{\fp}{\mathfrak{p}}

\newcommand{\fX}{\mathfrak{X}}

\newcommand{\V}{\mathcal{V}}
\newcommand{\cA}{\mathcal{A}}

\newcommand{\cD}{\mathcal{D}}
\newcommand{\ti}{\mathtt{i}}

\renewcommand{\emptyset}{\varnothing}
\newcommand{\tr}{\mathrm{tr}}
\newcommand{\CP}{\CC\mathrm{P}}
\renewcommand{\O}{\mathcal{O}}
\newcommand{\Lat}{\mathrm{Lat}}
\newcommand{\Gr}{\mathrm{Gr}}
\newcommand{\Fl}{\mathrm{Fl}}
\newcommand{\tI}{{\tt I}}
\newcommand{\tti}{{\tt i}}
\newcommand{\Omalg}{\Omega_\alg}
\newcommand{\Knarr}{{\mathrm{Knarr}}}
\newcommand{\alg}{{\mathrm{alg}}}
\newcommand{\rat}{{\mathrm{rat}}}
\newcommand{\diff}{{\mathrm{diff}}}
\newcommand{\fin}{{\mathrm{fin}}}
\newcommand{\weak}{{\mathrm{weak}}}
\newcommand{\Sut}{{\mathrm{S}\ut\,}}
\newcommand{\zpz}{z\partial_z}
\newcommand{\tl}{{\mathrm{tls}}}
\newcommand{\one}{\mathbf{1}}
\newcommand{\spa}{\mathrm{span}}
\newcommand{\Aut}{\mathrm{Aut}}
\newcommand{\Spe}{\mathrm{Spe}}
\newcommand{\Res}{\mathrm{Res}}
\newcommand{\PG}{\mathrm{PG}}
\newcommand{\Type}{\mathrm{type}}
\newcommand{\card}{\mathrm{card}}
\newcommand{\Sym}{\mathrm{Sym}}
\newcommand{\Cham}{\mathrm{Cham}}
\newcommand{\proj}{\mathrm{proj}}
\newcommand{\op}{\mathrm{op}}
\newcommand{\cC}{\mathcal{C}}
\newcommand{\cG}{\mathcal{G}}
\newcommand{\cS}{\mathcal{S}}
\newcommand{\kernel}{\mathrm{ker}}
\newcommand{\leqBruh}{\preceq}
\newcommand{\ut}{\scalebox{0.8}{\setlength{\unitlength}{1pt}%
\begin{picture}(10,10)%
\put(0,10){\line(1,0){10}}%
\put(10,10){\line(0,-1){10}}%
\put(0,10){\line(1,-1){10}}%
\end{picture}}}
\newcommand{\simplex}{{\scalebox{1.3}{\color{gray}\blacktriangle}}}
\newcommand{\fc}{{\scalebox{1.3}{\color{gray}\blacktriangleleft}}}
\begin{document}

\title{\bf Loop Groups and Twin Buildings}
\author{Linus Kramer\thanks{Supported by a Heisenberg fellowship
by the Deutsche Forschungsgemeinschaft}\\
\small Mathematisches Institut,
Universit\"at W\"urzburg,
Am Hubland,
D--97074 W\"urzburg,
Germany \\
\small email: {\tt kramer@mathematik.uni-wuerzburg.de}}
\date{}
\maketitle
\CompileMatrices
\begin{center}
\emph{Dedicated to John Stallings on the occasion of his 65th birthday.}
\end{center}
\begin{abstract}
In these notes we describe some buildings related to complex
Kac-Moody groups. First we describe the spherical building of
$\SL_n(\CC)$ (i.e. the projective geometry
$\PG(\CC^n)$) and its Veronese representation. Next we recall
the construction of the affine building associated to a
discrete valuation on the rational function field $\CC(z)$.
Then we describe the same building in terms of complex Laurent
polynomials, and introduce the Veronese representation, which
is an equivariant embedding of the building into an affine
Kac-Moody algebra. Next, we introduce topological twin buildings.
These buildings can be used for a proof
--- which is a variant of the proof by
Quillen and Mitchell \cite{Mitchell} ---
of Bott periodicity which uses only topological
geometry. At the end we indicate very briefly that the whole
process works also for affine real almost split Kac-Moody groups.
\end{abstract}


\section{Introduction}

We briefly recall the definition of Coxeter complexes and
buildings; for more details, we refer to the books by
Brown \cite{Brown}, Ronan \cite{Ronan}, Scharlau \cite{Scharlau},
and Tits \cite{Tits}. For our purposes, a \emph{simplicial complex}
is poset $(\Delta,\leq)$ whose elements are called simplices, with
the following two properties: any two simplices $X,Y\in\Delta$
have a unique infimum $X\sqcap Y$, and for any $X\in\Delta$,
the poset $\Delta_{\leq X}=\{Y\in\Delta|\ Y\leq X\}$ is
order-isomorphic to the power set $(2^F,\subseteq)$ of some finite set $F$;
the cardinality of this set $F$ is called the \emph{rank}
of the simplex~$X$. Note that the rank of a simplex differs
by one from the dimension of its geometric realization;
a $k$-simplex has rank $k+1$. The rank of a simplicial complex
is the maximum of the ranks of its simplices.

\begin{Num}
\textbf{Coxeter complexes}
\label{CoxeterComplex}
Let $I$ be a finite set with $r$ elements, and let
$(m_{ij})$ be a symmetric matrix indexed by $I\times I$, with entries
in $\mathbb{N}\cup\{\infty\}$, subject to the following two conditions:
$m_{ij}\geq 2$ for all $i\neq j$, and $m_{ii}=1$ for all $i$. 
Such a Coxeter matrix is determined by its Coxeter graph;
the vertices of this graph are the elements of $I$, and two
vertices $i,j$ are joined by $m_{ij}-2$ edges, or by one edge
labeled $m_{ij}$.
The corresponding \emph{Coxeter system} $(W,S)$ is the group $W$
with generating set $S=\{s_i|\ i\in I\}$ and relators
$(s_is_j)^{m_{ij}}$. Associated to such a Coxeter system is a
simplicial complex, the \emph{Coxeter complex} $\Sigma=\Sigma(W,S)$
which is defined as follows. For $J\subseteq I$ let
$W_J$ denote the subgroup generated by the elements $s_j$, for $j\in J$.
The simplices of $\Sigma$ are the cosets $wW_J\in W/W_J$, where $J$ runs
over all subsets of $I$; the partial ordering
is the reversed inclusion, i.e. 
\[
gW_J\leq hW_{J'}\text{ if and only if }
gW_J\supseteq hW_{J'}.
\]
The group $W$ acts regularly on the maximal simplices of
this simplicial complex (by
left translations). The \emph{type} of a coset $wW_J$ is
\[
\Type(wW_J)=I\setminus J;
\]
a subset $J\subseteq I$ is called \emph{spherical}
if $W_J$ is finite; if $I$ itself is spherical, then $\Sigma$ is
called spherical.
The geometric realization of a spherical Coxeter
complex of rank $k$ is a triangulated $k-1$-sphere.
We define a $W$-invariant double-coset valued distance function
$\delta$ on $\Sigma$ as follows:
\begin{align*}
\delta:\Sigma\times\Sigma\too&
\textstyle\bigcup\{W_J\backslash W/W_K|\ J,K\subseteq I\}\\
\delta(u W_J,vW_K) =&W_Ju^{-1}vW_K.
\end{align*}
\end{Num}
Coxeter groups have nice geometric properties; in particular, the
word problem can be solved. Let $A$ be an abelian group, and let
$a:S\too A$ be a function. We require that $a(s_i)=a(s_j)$ holds
whenever $m_{ij}$ is finite and odd. Then there is a well-defined
extension $a:W\too A$, which is defined as
$a(s_{i_1}\cdots s_{ir})=a(s_{i_1})+\cdots+a(s_{i_r})$, for
a \emph{reduced} (minimal) expression $w=s_{i_1}\cdots s_{i_r}$, the
\emph{$a$-length}. In the special case $A=\ZZ$, with $a(s_i)=1$ for
all $i$, we obtain the usual \emph{length function}
\[
\ell:W\too\ZZ .
\]
In general, the set of generators $S$ is not uniquely determined by
the abstract group $W$, so it is important to consider the pair $(W,S)$;
the question to which extent $S$ is determined by the group $W$ 
is treated in M\"uhlherr \cite{MueCox1} and
Brady, McCammond, M\"uhlherr \& Neumann \cite{MueCox2},
cp.~also Charney \& Davis \cite{CD} for related results.
\begin{Num}
\textbf{Buildings}
Let $\Delta\neq\emptyset$ be a simplicial complex, and let
$\Sigma=\Sigma(W,S)$
be a Coxeter complex. A simplicial injection $\phi:\Sigma\too\Delta$
is called a \emph{chart}, and its image $A=\phi(\Sigma)$ is
called an \emph{apartment}. The complex $\Delta$ is called a 
\emph{building} (of type $(m_{ij})$ and rank $r$)
if there exists a collection $\cA$ of apartments with the
following properties.

\smallskip\textbf{Bld}$_1$
For any two simplices $X,Y\in\Delta$, there exists an apartment
$A\in\cA$ with $X,Y\in A$.

\smallskip\textbf{Bld}$_2$
Given two charts $\phi_i:\Sigma\rInto\Delta$, for
$i=1,2$, there exists an element $w\in W$ such that
$\phi_1\circ w(X)=\phi_2(X)$ holds for all
$X\in\phi^{-1}_2(\phi_1(\Sigma)\cap\phi_2(\Sigma))$.

\smallskip
A simplex of maximal rank $r$ is also called a \emph{chamber}; the
set of all chambers is $\Cham(\Delta)$.
For any subset $X\subseteq \Delta$, we let $\Cham(X)$ denote the set of
all chambers contained in $X$.
The building is \emph{thick} if every simplex of rank $r-1$
is contained in at least three chambers.
\textbf{All buildings in this paper will be thick.}
A \emph{spherical building} is a building with finite apartments.

It follows from the axioms that there is a well-defined double-coset
valued \emph{distance function}
\[
\delta:\Delta\times\Delta\too
\textstyle\bigcup\{W_J\backslash W/W_K|\ J,K\subseteq I\}
\]
whose restriction to any apartment is given by the function $\delta$
defined above in \ref{CoxeterComplex}.
The restriction of $\delta$ to $\Cham(\Delta)\times\Cham(\Delta)$
is Tits' more familiar $W$-valued distance function; buildings can
also be characterized by properties of the distance function $\delta$,
see Ronan's book \cite{Ronan}.
\end{Num}
Let $\simplex^{r-1}$ denote the standard $r-1$-simplex $(2^I,\subseteq)$.
The two axioms yield a simplicial surjection, the 
\emph{type function} (the 'accordion map')
\[
\Type:\Delta\too\simplex^{r-1},
\]
whose restriction to any apartment agrees with the type function
defined above. The type function is characterized (up to automorphisms
of $\simplex^{r-1}$) by the fact that its restriction to every
simplex of $\Delta$ is injective. A simplex $X\in\Delta$ is called
\emph{spherical} if $I\setminus\Type(X)$ is spherical.

For non-spherical buildings the apartment
system $\cA$ is in general not unique
(there is always a unique \emph{maximal} apartment system),
but the isomorphism
type of the apartments and the Coxeter system $(W,S)$ is uniquely
determined by the simplicial complex $\Delta$.

\begin{Num}
\textbf{Automorphisms}
The automorphism group $\Aut(\Delta)$ consists
of all simplicial automorphisms of $\Delta$; it has a normal
subgroup $\Spe(\Delta)$ consisting of all type-preserving
automorphisms. An \emph{action} of a group $G$ on $\Delta$
is a homomorphism $G\too\Aut(\Delta)$. We say that $G$ acts
\emph{transitively} on $\Delta$ if it acts transitively on the
set $\Cham(\Delta)$ of chambers.

\medskip\noindent\textbf{Strongly transitive actions and $BN$-pairs}
A group $G$ is said to act \emph{strongly transitively} on $\Delta$
(with respect to $\mathcal A$) if $G$ acts as a group of special
automorphisms on $\Delta$ such that

\smallskip
\textbf{STA}$_1$ $G$ acts transitively on the set $\mathcal A$
of apartments.

\smallskip
\textbf{STA}$_2$ If $A\in\cA$ is an apartment,
then the set-wise stabilizer $N$ of $A$ acts transitively on the 
chambers in $A$.

\smallskip
Let $C\in A$ be a chamber in an apartment $A$, let $B=G_C$ denote
the stabilizer of $C$, and $N$ the set-wise stabilizer of $A$. Then
$(B,N)$ is a so-called $BN$-pair for the group $G$; the Weyl
group $W=N/(N\cap B)$ acts regularly on $A$ and
is isomorphic to the Coxeter group
of $\Delta$. The stabilizers of the simplices contained in $C$
are called \emph{standard parabolic subgroups} of $G$.
\end{Num}

\begin{Num}
\textbf{Panels and residues}
Let $X\in\Delta$ be a simplex of type $J$. The \emph{residue} of $X$
is the poset
$\Res(X)=\Delta_{\geq X}=\{Y\in\Delta|\ Y\geq X\}$.
The residue is order-isomorphic to the \emph{link} of $X$, and
thus can be identified with a subcomplex of $\Delta$
(although strictly speaking, $\Res(X)$ is not a subcomplex).
It is a basic but important fact that residues are again
buildings; the corresponding Coxeter complex is obtained
from the restricted Coxeter matrix $(m_{i.j})_{J\times J}$.
The \emph{type}
of the residue $\Res(X)$ is $I\setminus J$, and its rank
is $\card(I\setminus J)$. A residue $S_i=\Res(X)$ of 
rank 1 and type $\{i\}$ is called an $i$-\emph{panel}.
\end{Num}
The following observation is very simple, but useful.

\begin{Lem}
\label{TransLemma}
Suppose that $G$ acts as a group of special automorphisms on
a building $\Delta$ of rank $r\geq 2$; let $C\in\Delta$ be a chamber.
Then $G$ acts transitively on $\Delta$ if and only if the following
holds for all $r$ simplices $X_1,\ldots,X_r\leq C$ of rank $r-1$:

\smallskip
$\mathbf{\bullet}$ the stabilizer $G_{X_i}$ acts transitively on the
panel $\Res(X_i)$.
\qed
\end{Lem}
Tits pointed out that
transitive groups acting on buildings of higher rank can be represented
as amalgams.

\begin{Thm}
(Tits \cite{StAndrews} 2.3)\label{Amalgam}
Let $G$ be a group acting transitively and type-preservingly on an
irreducible building $\Delta$ of rank $r\geq 3$ (i.e. we assume that the
Coxeter diagram is connected and has at least 3 vertices).
Let $X$ be a chamber,
let $X_1,\ldots,X_r$, denote its subsimplices of rank $r-1$, and
$X_{\{ij\}}$, $i,j\in I$, $i\neq j$ its $\binom{r}{2}$ subsimplices of rank
$r-2$. Then the $\binom{r+1}{2}+1$ 
different $G$-stabilizers of these simplices with their natural inclusions
form a diagram (a simple 2-complex of groups,
cp.~Bridson \& Haefliger~\cite{BridHaef}  II.12 and
III.$\mathcal{C}$
--- the corresponding poset is the set of all subsets of
$I$ with at most two elements) whose limit is $G$.
\qed
\end{Thm}
For example, a building of rank 4 with
$I=\{1,2,3,4\}$ yields a complex of groups
as follows, if we put $G_{X_J}=G_J$ for short.
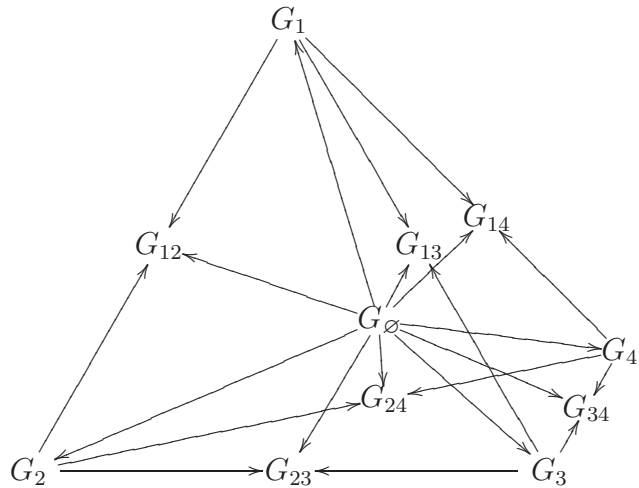
\[
\begin{xy} /r4cm/:
{\xypolygon3"B"{~*{\phantom{G_{12}}}~>{}}}
,{\xypolygon3"C"{~:{(0.5,0):}~={30}~*{\phantom{G_1}}~>{}}}
,{(1.3,0)*{G_\emptyset}="O"}
,{(2.1,-0.1)*{G_4}="B4"}
,{"B1"*{G_1}}
,{"B2"*{G_2}}
,{"B3"*{G_3}}
,{"C1"*{G_{13}}}
,{"C2"*{G_{12}}}
,{"C3"*{G_{23}}}
,{"B1";"B4"**@{},?<>(0.6)*{G_{14}}="C4"}
,{"B2";"B4"**@{},?<>(0.6)*{G_{24}}="C5"}
,{"B3";"B4"**@{},?<>(0.5)*{G_{34}}="C6"}
,{\ar @{->}"B1";"C1"}
,{\ar @{->}"B1";"C2"}
,{\ar @{->}"B2";"C2"}
,{\ar @{->}"B2";"C3"}
,{\ar @{->}"B3";"C1"}
,{\ar @{->}"B3";"C3"}
,{\ar @{->}"O";"C1"}
,{\ar @{->}"O";"C2"}
,{\ar @{->}"O";"C3"}
,{\ar @{->}"O";"C4"}
,{\ar @{->}"O";"C5"}
,{\ar @{->}"O";"C6"}
,{\ar @{->}"O";"B1"}
,{\ar @{->}"O";"B2"}
,{\ar @{->}"O";"B3"}
,{\ar @{->}"O";"B4"}
,{\ar @{->}"B1";"C4"}
,{\ar @{->}"B4";"C4"}
,{\ar @{->}"B2";"C5"}
,{\ar @{->}"B4";"C5"}
,{\ar @{->}"B3";"C6"}
,{\ar @{->}"B4";"C6"}
\end{xy}
\]
In general, the geometric realization of this 2-complex of groups
is the cone over the first barycentric subdivision of the 1-skeleton
of an $r-1$-simplex.

For group amalgamations related to twin buildings (cp.~Section
\ref{TwinBuildings}) see Abramenko \& M\"uhl\-herr~\cite{AM}.

\section{The spherical case: projective space and $\SL_n(\CC)$}
\label{SphericalCase}
In this section we describe the spherical building
obtained from $n-1$-dimensional complex projective space,
some related groups, and the Veronese representation.
Given a ring $R$, we let $R(k)$ denote the matrix algebra of
all $k\times k$-matrices with entries in $R$.

\subsection*{The spherical building $\Delta(\CC^n)$}

The \emph{Grassmannian} of $k$-spaces in $\CC^n$ is the 
complex projective variety
\[
\Gr_k(\CC^n)=\{U\leq \CC^n|\ \dim(U)=k\}.
\]
A \emph{partial flag} in $\CC^n$ is a nested sequence of subspaces
$0<U_{j_1}<U_{j_2}<\cdots<U_{j_r}<\CC^n$, with $\dim(U_{j_\nu})=j_\nu$.
Such a partial flag can be viewed as a map 
\[
\{1,\ldots,n-1\}\supseteq J\too^{U}
\Gr_1(\CC^n)\cup\cdots\cup\Gr_{n-1}(\CC^n),
\]
with $\dim(U_j)=j$ and $U_j\leq U_k$ for $j\leq k$. There is a
natural order '$\leq$' on the collection $\Delta(\CC^n)$ of all
partial flags: if $U$ and $U'$ are partial flags with
$J\subseteq J'\subseteq \{1,\ldots,n-1\}$, then $U\leq U'$
if and and only if $U'|_J=U$.
The resulting poset is a simplicial complex of rank $n-1$,
the spherical building
\[
\Delta(\CC^n)=(\Delta(\CC^n),\leq).
\]
Of course, $\Delta(\CC^n)$ is precisely the same as the $n-1$-dimensional
projective geometry $\PG(\CC^n)$ in a different guise.

Given an ordered basis $(v_1,\ldots,v_n)$ of $\CC^n$, we may consider
the maximal flag $U=U^{\CC}(v_1,\ldots,\linebreak v_n)$, where
$U_i=\spa_\CC\{v_1,\ldots,v_i\}$, and the apartment
$A^\CC\{v_1,\ldots,v_n\}$
consisting of all flags obtained as partial
flags from the $n!$ distinct maximal flags (chambers)
$U^\CC(v_{\pi(1)},\ldots,v_{\pi(n)})$, where $\pi\in\Sym(n)$.
It is not difficult to check that $\Delta(\CC^n)$ together with
this collection of apartments is a spherical building.
As a simplicial complex, 
$A^\CC\{v_1,\ldots,v_n\}$ is a triangulation of the sphere $\SS^{n-2}$.
In fact, let $\simplex^{n-1}$ denote the standard $n-1$-simplex,
and $\mathrm{Bd}(\simplex^{n-1})$ its boundary; 
then $A^\CC\{v_1,\ldots,v_n\}$ is simplicially isomorphic to
the first barycentric subdivision
$\Sigma=\mathrm{Sd}(\mathrm{Bd}(\simplex^{n-1}))$.
The corresponding Coxeter group is
the symmetric group $W=\Sym(n)$, given by the presentation
\[
W=\langle s_1,\ldots,s_{n-1}|\ (s_is_j)^{m_{ij}}=1,\ 1\leq i,j\leq n-1
\rangle,
\text{ where }
m_{ij}=\begin{cases} 1 &\text{ if }i=j\\
3 &\text{ if }|i-j|=1\\
2 &\text{ else.}
\end{cases}
\]
The involution $s_i$ is the transposition $(i,i+1)$; the Coxeter diagram
is
\[
A_{n-1}:\qquad
\begin{xy}
*@{*} \PATH ~={**\dir{-}}
'(5,0)*@{*}
'(10,0)*@{*} 
'(15,0)*@{*} 
'(20,0)*@{*} 
'(25,0)*@{*}
'(30,0)*@{*}
'(35,0)*@{*} 
'(40,0)*@{*} 
'(45,0)*@{*} 
'(50,0)*@{*}
\end{xy}
\]
($n-1$ nodes).
The type of a flag
\[
U:J\too\Gr_1(\CC^n)\cup\cdots\cup\Gr_{n-1}(\CC^n)
\]
is $\Type(U)=J$, the
set of the dimensions of the subspaces occurring in $U$.
In $\Delta(\CC^n)$, a panel $S_i$ is determined by a flag of the form 
\[
U_1<U_2<\cdots<U_{i-1}<U_{i+1}<\cdots U_{n-1},
\]
and there is a natural bijection
$S_i\too \Gr_1(U_{i+1}/U_{i-1})\cong\CC\mathrm{P}^1=\CC\cup\{\infty\}$
onto a projective line.

\subsection*{Parabolics in $\SL_n(\CC)$}

There are natural actions of the groups $\GL_n(\CC)$ and $\SL_n(\CC)$ on
the building $\Delta(\CC^n)$ (by type preserving automorphisms), and
it is not difficult to see that these actions are strongly transitive.
The $\SL_n(\CC)$-stabilizer of a partial flag is thus a parabolic
subgroup.
The \emph{maximal parabolics} are the stabilizers of the flags
of rank $1$, i.e. of subspaces $0<V<\CC^n$.
Such a maximal parabolic is conjugate to one of the
\emph{standard maximal parabolics}
\[
P^i=\left.\left\{\begin{pmatrix} A & B \\ 0 & C \end{pmatrix}\right|\ 
A\in\CC(i),\ C\in\CC(n-i),\ B\in\CC^{i\times (n-i)},\ \det(A)\det(C)=1
\right\},
\]
and
\[
\Gr_i(\CC^n)\cong\SL_n(\CC)/P^i.
\]
The minimal parabolics or \emph{Borel subgroups} are the
stabilizers of maximal flags, i.e. chambers in $\Delta(\CC^n)$.
The Borel subgroups are conjugate to the \emph{standard Borel subgroup}
\[
B=P^1\cap \cdots \cap P^{n-1}=\Sut_n(\CC)=
\left.\left\{
\begin{pmatrix} a_1 && * \\ & \ddots & \\ 0 && a_n \end{pmatrix}
\right|\ 
a_1a_2\cdots a_n=1\right\}
\]
consisting of all unimodular upper triangular matrices.
The corresponding homogeneous space is the \emph{complex flag variety}
\[
\Fl(\CC^n)=\Cham(\Delta(\CC^n))=
\{(U_1,\ldots,U_{n-1})|\ 0<U_1<\cdots<U_{n-1}<\CC^n\}
\cong\SL_n(\CC)/B.
\]

\subsection*{The anisotropic real structure of $\SL_n(\CC)$}
\label{FiniteCartan}
We define a semi-linear involution $*$ on $\CC(n)$ by putting
$X^*=\bar X^T$ (conjugate transpose). Recall that the Lie algebra
$\fsl_n\CC$ consists of all traceless $n\times n$-matrices.
Let 
\[
\fsu(n)=\{X\in\fsl_n\CC|\ X+X^*=0\}\quad\text{and}\quad
\fp_n=\{X\in\fsl_n\CC|\ X=X^*\}.
\]
The decomposition $\fsl_n(\CC)=\fsu(n)\oplus\fp_n$ is called the
\emph{Cartan decomposition} of the Lie algebra $\fsl_n(\CC)$.
On the group level, let 
\[
\SU(n)=\{g\in\SL_n(\CC)|\ g^*=g^{-1}\}
\]
denote the group of fixed elements of the involution
$(g\mapstoo g^{-*})\in\Aut_\RR(\SL_n(\CC))$.
Both the involution $X\mapstoo -X^*$
in the Lie algebra and the involution $g\mapstoo g^{-*}$ in the
group are called \emph{Cartan involutions}.
(In terms of algebraic groups, we have defined an
\emph{almost simple $\RR$-group scheme $\underline G$} such that
$\underline G(\RR)=\SU(n)$ and $\underline G(\CC)=\SL_n(\CC)$.
The Galois group of $\CC/\RR$ acts on the group 
$\SL_n(\CC)$ of $\CC$-points
of $\underline G$, and $\SU(n)$ is the group of fixed elements.
The group scheme $\underline G$ is $\RR$-\emph{anisotropic}:
no parabolic of $\underline G$ is defined over $\RR$.
For semi-simple $\RR$-algebraic groups, 'anisotropic' is the same as
'compact'.)

More geometrically, the group $\SU(n)$ can be described as follows.
Let $\langle x,y\rangle=\bar x_1y_1+\cdots+\bar x_ny_n$ denote the
standard hermitian form on $\CC^n$. This form induces a map
$\Gr_k(\CC^n)\too\Gr_{n-k}(\CC^n)$, $V\mapstoo V^\perp$, for all $k$,
which extends
in a natural way to an involution $\perp$ on $\Delta(\CC^n)$. Then
$\SU(n)$ is the centralizer of this involution,
\[
\SU(n)=\mathrm{Cen}_{\SL_n(\CC)}(\perp).
\]
Classically, the involution $\perp$ is called the
\emph{standard elliptic polarity} on the complex projective
geometry $\PG(\CC^n)$; the associated Riemannian symmetric space
$X=\SL_n(\CC)/\SU(n)$ is --- as a subset of $\Aut(\Delta)$ ---
precisely the space of all elliptic polarities
(i.e. the space of all positive definite hermitian forms on $\CC^n$).

Using Gram-Schmidt orthonormalization, one shows that
$\SU(n)$ acts transitively on the flag variety
$\Fl(\CC^n)$. Bernhard M\"uhlherr pointed out that there is a different
proof which uses Lemma \ref{TransLemma}, and which carries
over to the Kac-Moody case described later.
\begin{Lem}
\label{SU(n)IsTransitive}
The group $\SU(n)$ acts transitively on $\Fl(\CC^n)$.

\proof
Let $U=U^\CC(e_1,\ldots,e_n)=(U_1,\ldots,U_{n-1})$ denote the maximal
flag arising from the standard basis of $\CC^n$. The $\SU(n)$-stabilizer
of a panel $S_i=\Res(U_1,\ldots,U_{i-1},U_{i+1},\ldots,U_{n-1})$
containing $U$ is isomorphic to
$\SU(2)\cdot T^{n-2}$; it induces the transitive
group $\SO(3)=\SU(2)\cdot T^{n-2}/T^{n-2}$ on the panel
(we let $T^k=\U(1)\times\cdots\times\U(1)$ denote the compact
torus of rank $k$; in terms of algebraic groups, $T^k$ is the group 
$\underline T(\RR)$ of
$\RR$-points of an anisotropic $\RR$-torus $\underline T$ of rank $k$).
The assertion follows with Lemma \ref{TransLemma}.
\qed
\end{Lem}
For $n=4$, we have by Theorem \ref{Amalgam}
the following complex of groups which
represents $\SU(4)$ as an amalgam (we indicate the nonzero entries
of a matrix by $*$, and all matrices are assumed to be unimodular
and unitary).
\[
\begin{xy} /r3cm/:
*{\scalebox{0.4}{\begin{pmatrix}*&&&\\&*&&\\&&*&\\&&&*\end{pmatrix}}}="O"
,{\xypolygon3"B"{~*{\phantom{%
\scalebox{0.4}{\begin{pmatrix}*&&&\\&*&&\\&&*&\\&&&*\end{pmatrix}}}}~>{}}}
,{\xypolygon3"C"{~:{(0.5,0):}~={30}~*{\phantom{%
\scalebox{0.4}{\begin{pmatrix}*&&&\\&*&&\\&&*&\\&&&*\end{pmatrix}}}}~>{}}}
,"B3"*{\scalebox{0.4}{\begin{pmatrix}*&*&&\\{*}&*&&\\&&*&\\&&&*\end{pmatrix}}}
,"B2"*{\scalebox{0.4}{\begin{pmatrix}*&&&\\&*&*&\\&*&*&\\&&&*\end{pmatrix}}}
,"B1"*{\scalebox{0.4}{\begin{pmatrix}*&&&\\&*&&\\&&*&*\\&&*&*\end{pmatrix}}}
,"C3"*{\scalebox{0.4}{\begin{pmatrix}*&*&*&\\{*}&*&*&\\{*}&*&*&\\&&&*\end{pmatrix}}}
,"C1"*{\scalebox{0.4}{\begin{pmatrix}*&*&&\\{*}&*&&\\&&*&*\\&&*&*\end{pmatrix}}}
,"C2"*{\scalebox{0.4}{\begin{pmatrix}*&&&\\&*&*&*\\&*&*&*\\&*&*&*\end{pmatrix}}}
,{\ar @{->}"B1";"C1"}
,{\ar @{->}"B1";"C2"}
,{\ar @{->}"O";"C1"}
,{\ar @{->}"B2";"C2"}
,{\ar @{->}"B2";"C3"}
,{\ar @{->}"O";"C2"}
,{\ar @{->}"B3";"C3"}
,{\ar @{->}"B3";"C1"}
,{\ar @{->}"O";"C3"}
,{\ar @{->}"O";"B1"}
,{\ar @{->}"O";"B2"}
,{\ar @{->}"O";"B3"}
\end{xy}
\]
Here
\begin{align*}
\scalebox{0.4}{\begin{pmatrix}*&*&&\\{*}&*&&\\&&*&\\&&&*\end{pmatrix}}\cong
\scalebox{0.4}{\begin{pmatrix}*&&&\\&*&*&\\&*&*&\\&&&*\end{pmatrix}}\cong
\scalebox{0.4}{\begin{pmatrix}*&&&\\&*&&\\&&*&*\\&&*&*\end{pmatrix}}\cong
\SU(2)\cdot T^2\\
\scalebox{0.4}{\begin{pmatrix}*&*&*&\\{*}&*&*&\\{*}&*&*&\\&&&*\end{pmatrix}}
\cong
\scalebox{0.4}{\begin{pmatrix}*&&&\\&*&*&*\\&*&*&*\\&*&*&*\end{pmatrix}}
\cong\SU(3)\cdot T^1\\
\scalebox{0.4}{\begin{pmatrix}*&*&&\\{*}&*&&\\&&*&*\\&&*&*\end{pmatrix}}
\cong\SU(2)\cdot\SU(2)\cdot T^1.
\end{align*}
In group theoretic terms, we thus have
\begin{align*}
\SL_n(\CC)=\SU(n)B & \qquad \SU(n)\cap B=T^{n-1} \\
\SL_n(\CC)=\SU(n)P^i & \qquad \SU(n)\cap P^i\cong
\mathrm{S}(\U(i)\times\U(n-i)).
\end{align*}
In particular, 
\[
\Fl(\CC^n)\cong\SU(n)/T^{n-1}\quad\text{and}\quad
\Gr_k(\CC^n)\cong\SU(n)/\mathrm{S}(\U(k)\cdot\U(n-k)).
\]

\subsection*{Knarr's construction}
Let $|\Delta(\CC^n)|$ denote the \emph{geometric realization} of the simplicial
complex $\Delta(\CC^n)$. There are various ways to topologize this set.
One possibility is the \emph{weak topology} determined by the simplices:
by definition, a subset $A\subseteq|\Delta(\CC^n)|$
is closed in the weak topology if and only if its intersection
with every simplex is closed. We denote the resulting topological
space by $|\Delta(\CC^n)|_\weak$; there are other nice
topologies, all of which yield the same weak homotopy type for
the space $|\Delta(\CC^n)|$, cp.~Bridson \& Haefliger
\cite{BridHaef}~I.7.
The Solomon-Tits Theorem asserts
in our situation that there is a homotopy equivalence
\[
|\Delta(\CC^n)|_\weak\simeq\SS^{n-2}\wedge X_+,
\]
where $X$ is a discrete set of cardinality $2^{\aleph_0}$ and $X_+$ its 
one-point compactification.
(In general, the Solomon-Tits Theorem says that for a spherical building 
$\Delta$ of rank $r$,
\[
|\Delta|_\weak\simeq\SS^{r-1}\wedge X_+\simeq\bigvee_{\card(X)}\SS^{r-1}
\]
where $X$ is a discrete space whose cardinality is
$\card(X)=\card\{A\in\cA|\ C\in A\}$ for
some fixed chamber $C$.
The action of the automorphism group $\SL_n(\CC)$ on the top-dimensional
homology group $H_{r-1}(|\Delta|_\weak)$ of this complex is called the
\emph{Steinberg representation.} See Ronan \cite{Ronan} App.~4 for more
details and further references.)

However, this construction neglects the natural topology of
$\Delta(\CC^n)$. Consider the following construction.
Fix an $(n-2)$-simplex
$\simplex^{n-2}$ and label its vertices as $1,2,\ldots,n-1$.
There is a natural surjection
\[
\Fl(\CC^n)\times|\simplex^{n-2}|\too|\Delta(\CC^n)|
\]
which maps $\{U\}\times|\simplex^{n-2}|$ to the geometric
realization of the simplex of $\Delta(\CC^n)$
spanned by the vertices $U_1,\ldots,U_{n-1}$ of the given flag $U$,
in such a way that the
$i$th vertex of $\{U\}\times\simplex^{n-2}$ is identified with $U_i$.
There is a natural compact
topology on $\Fl(\CC^n)\times|\simplex^{n-2}|$,
and we endow $|\Delta(\CC^n)|$ with the quotient topology.
We denote the resulting space by $|\Delta(\CC^n)|_\Knarr$
(because Knarr --- inspired by Mitchell \cite{Mitchell} ---
introduced it first in \cite{Knarr} for compact buildings of rank 2).
It can be shown that there is a homeomorphism
\[
|\Delta(\CC^n)|_\Knarr\cong\SS^{n^2-2},
\]
see Knarr \emph{loc.cit.}
We will prove this in the next section, using the Veronese
representation of $\Delta(\CC^n)$; a more general result is stated
in Section \ref{TopologicalGeometry}.

The Knarr construction works for general topological
buildings. If the topology on the
chamber set of a spherical building $\Delta$ of rank $r$
satisfies certain natural
conditions (e.g.~the inclusions between its Schubert varieties should be
cofibrations), then $|\Delta|_\Knarr\simeq\SS^{r-1}\wedge O_+$,
where $O$ is the set of all chambers opposite to a fixed chamber,
and $O_+$ its one-point compactification. For the special case
of a discrete topology on the chamber set, this is precisely the
Solomon-Tits Theorem. There is a well-developed theory of
compact spherical buildings (the case of rank 2 is worked out in
Kramer \cite{KraDiss}, and the results proved there extend immediately
to the case of higher rank); the result for the homotopy
type of $|\Delta|_\Knarr$ can be proved in much greater generality,
see \cite{KraDiss} Sec.~3.3.
If the building is spherical, irreducible, compact,
connected, and of rank at least 3
(see Section \ref{TopologicalGeometry} for definitions),
then by the results of Burns \& Spatzier \cite{BurnsSpatzier} and 
Grundh\"ofer, Knarr \& Kramer
\cite{Flaghom1},
\cite{Flaghom2},
\cite{Flaghom3},
the space
$|\Delta|_\Knarr$ can be identified with the 
\emph{visual boundary} $X(\infty)$, cp.~Bridson \& Haefliger
\cite{BridHaef}~II.8.,
of a Riemannian symmetric space $X$; for $\Delta(\CC^n)$, the symmetric
space in question is $X=\SL_n(\CC)/\SU(n)$;
the same conclusion holds for buildings of rank 2, provided that
the automorphism group acts transitively on the flags,
see Grundh\"ofer, Knarr \& Kramer
\cite{Flaghom1},
\cite{Flaghom2},
\cite{Flaghom3}.

\subsection*{The Veronese representation of $\Delta(\CC^n)$}
We endow $\CC(n)$ with the positive definite hermitian form
\[
\langle X,Y\rangle=\tr(X^*Y).
\]
Consider the subspace $H(n)\leq\CC(n)$ consisting of
all hermitian matrices in $\CC(n)$.
Every matrix $X\in H(n)$ has a unique decomposition
\[
X=X^\tl+\textstyle\frac{\tr(X)}n\one.
\]
into a traceless hermitian matrix $X^\tl$
and a real multiple of the identity
matrix.
The adjoint action of $\SU(n)$ on $H(n)$ is the action by
conjugation,
\[
X\mapstoo gXg^*.
\]
We have an orthogonal $\SU(n)$-invariant
splitting $H(n)=\fp_n\oplus \one\RR$ (recall that
$\fp_n=\{X\in H(n)|\linebreak\tr(X)=0\}$).
Suppose that $X\in H(n)$ is a \emph{projector}, i.e. that
\[
X^2=X.
\]
If $X\neq0,\one$, then
the minimal polynomial of $X$ is $\mu_X(t)=t(t-1)$, and 
$\tr(X)=k$, for some $k\in\{1,\ldots,n-1\}$.
The kernel $V$ of $X$ is then an $n-k$-dimensional 
subspace of $\CC^n$. In this way, we obtain an $\SU(n)$-equivariant
1-1 correspondence between elements of $\Gr_{n-k}(\CC^n)$
and self-adjoint projectors with trace $k$ which is given
by the map $X\mapstoo \kernel(X)$.
The map $X\mapstoo X^\tl=X-\frac{\tr(X)}n\one$ is 
$\SU(n)$-equivariant; in this way, we obtain an embedding 
$\Phi:\Gr_k(\CC^n)\rInto\fp_n$
as follows. For $V\in\Gr_k(\CC^n)$ let $X_V$ denote the unique 
self-adjoint projector
with $\kernel(X_V)=V$, and put $\Phi(V)=(X_V)^\tl$.
The elliptic polarity $\perp$ is built-in: the other eigenspace
of $X_V$ is the image $V^\perp$ of $V$ under the
elliptic polarity.
Even better, the incidence can be seen in $\fp_n$:
two self-adjoint operators $\Phi(V),\Phi(W)\in\fp_n$ representing subspaces
$V\in\Gr_i(\CC^n)$ and $W\in\Gr_j(\CC^n)$ are incident if and only
if the euclidean distance $|\Phi(V)-\Phi(W)|$ attains the minimum
possible value
$d_{ij}=\mathrm{dist}(\Phi(\Gr_i(\CC^n),\Phi(\Gr_j(\CC^n))$.
If $U$ is a flag in $\Delta(\CC^n)$,
and if $p\in\simplex^{n-2}$ has barycentric coordinates
$(p_1,\cdots,p_{n-1})$, then
we map $(U,p)\in\Fl(\CC^n)\times\simplex^{n-2}$ to the 
hermitian operator
\[
\Phi(U,p)=\sum_{i=1}^{n-1} p_i\Phi(U_i)\in\fp_n.
\]
In this way we obtain an $\SU(n)$-equivariant injection
\[
\Phi:|\Delta(\CC^n)|\rInto\fp_n.
\]
We call this the \emph{Veronese representation} of the building
$\Delta(\CC^n)$. The image of the flag space $\Fl(\CC^n)$
in $\fp_n$ is an \emph{isoparametric submanifold}
(we identify a chamber with the barycenter of its geometric
realization); the images of the partial flag varieties are parallel focal
submanifolds in this isoparametric foliation, cp.~Palais \& Terng
\cite{PalaisTerng}, Thorbergsson \cite{Thorbg} and \cite{ThorbgHandbook},
Knarr \& Kramer \cite{KK}.
(The corresponding construction for the real projective
geometry $\PG(\RR^3)$ leads to the classical 
\emph{Veronese embedding} of $\RR\mathrm{P}^2\rInto\SS^4$, whence the name.)

The following variation of the map $\Phi$ is also useful.
For a nonzero matrix $X\in\CC(n)$, put $\widehat X=|X|^{-1}X$,
where $|.|$ denotes the euclidean norm,
and consider the map $\widehat\Phi:(U,p)\mapstoo\widehat{\Phi(U,p)}$.
Then is is not difficult to see that
\[
\widehat\Phi(|\Delta(\CC^n)|)=\SS^{n^2-2}\subset\fp_n\cong\RR^{n^2-1};
\]
the map is injective, since we can recover a flag $(U_{i_1},\ldots,U_{i_r})$
from its image $X\in\fp_n$ as follows: the hermitian matrix $X$ has
eigenvalues $\lambda_1<\lambda_2<\ldots<\lambda_r$, and
$U_{i_k}=\kernel(X-\lambda_1\one)
\oplus\cdots\oplus\kernel(X-\lambda_k\one)$.
The surjectivity follows from
the fact that every hermitian matrix can be diagonalized under the
$\SU(n)$-action, because the image of the Veronese representation
contains certainly all diagonal traceless matrices of norm 1; these
are precisely the images of the simplices in the
apartment $A^\CC\{e_1,\ldots,e_n\}$.
Since $\Phi$ is continuous on $|\Delta(\CC^n)|_\Knarr$,
we obtain in particular the claimed homeomorphism
\[
|\Delta(\CC^n)|_\Knarr\cong\SS^{n^2-2}.
\]

\section{The affine building of $\CC(z)$}

In this section we describe the affine building associated to the
discrete valuation on the rational function field $\CC(z)$.
This building is discussed in considerably more detail in the books by
Brown \cite{Brown} and Ronan \cite{Ronan}.

\subsection*{Lattices in $\CC(z)^n$}
We let $\LL=\CC(z)$ denote the field of fractions of the polynomial
ring $\CC[z]$. Thus, $\LL$ is the field
of rational functions on the complex projective line
$\CP^1=\CC\cup\{\infty\}$. For $c\in\CP^1$ we let
$\mathcal{O}_c\leq\LL$ denote the subring of all rational functions
which don't have a pole in $c$, and $m_c=\{f\in\O_c|\ f(c)=0\}$
the maximal ideal of $\O_c$. Evaluation at $c$ yields a map
$ev_c:\O_c\too\CC$, $f\mapstoo f(c)$ with kernel $m_c$, and we obtain
exact sequences
\[
0\too m_c\too \O_c\too^{ev_c}\CC\too0
\]
and
\[
1\too\SL_n(m_c)\too\SL_n(\O_c)\too^{ev_c}\SL_n(\CC)\too1
\]
(the group $\SL_n(m_c)$ is defined to be the kernel of the
evaluation map $ev_c$).
We may view the elements of $\SL_n(\LL)$ as
\emph{rational maps} from $\CP^1$ into $\SL_n(\CC)$.
Note also that every element $q\neq 0$
of $\LL$ can be expressed in the form
\[
\textstyle q=(z-c)^k\frac{f}{g},
\]
with $k\in\ZZ$, $f,g\in\CC[z]$, and $f(c)\neq0\neq g(c)$.
We put $\nu_c(q)=k$, and $\nu_c(0)=\infty$. The map 
\[
\nu_c:\LL\too\ZZ\cup\{\infty\}
\]
is a \emph{discrete valuation} on $\LL$
(with some modifications for $c=\infty$).
Note that $\O_c^\times=\{q\in\O_c|\ \nu_c(q)=0\}$.
There is nothing special about the choice of $c\in\CP^1$, and we put $c=0$
for the remainder of this section.

The group $\SL_n(\LL)$ acts on the projective geometry
$\PG(\LL^n)$ in very much the same way as $\SL_n(\CC)$ on $\PG(\CC^n)$,
and we could consider the spherical building $\Delta(\LL^n)$.
But now we introduce a different geometry for this
group, the \emph{affine building} $\Delta(\LL^n,\O_0)$.
Given an $\LL$-basis $v_1,\ldots,v_n$ of $\LL^n$, we have the
free $\O_0$-module 
\[
M=\spa_{\O_0}\{v_1,\ldots,v_n\}=v_1\O_0+\cdots+v_n\O_0
\]
of rank $n$ generated by
these basis vectors. We call $M$ an $\O_0$-\emph{lattice}, and we let
$\Lat_n(\LL,\O_0)$ denote the collection of
all such lattices.
(The following simple observation is useful. If
$M\subseteq\LL^n$ is a free $\O_0$-module of rank $k$, with $\O_0$-basis
$\{v_1,\ldots,v_k\}$, then $\{v_1,\ldots,v_k\}$ is linearly independent
over $\LL$, because $\LL$ is the field of fractions of $\O_0$. Thus,
the $\O_0$-lattices are precisely the free $\O_0$-modules of rank $n$
in $\LL^n$.)
Evidently, the group $\GL_n(\LL)$ acts transitively
on $\Lat_n(\LL,\O_0)$; the $\GL_n(\LL)$-stabilizer of the $\O_0$-module
$M_0$ spanned by the canonical basis $e_1,\ldots,e_n$ of $\LL^n$ is
the group $\GL_n(\O_0)$, and therefore
\[
\Lat_n(\LL,\O_0)\cong\GL_n(\LL)/\GL_n(\O_0).
\]
We call two lattices $M,M'\in\Lat_n(\LL,\O_0)$ \emph{projectively
equivalent} if $M=qM'$ for some $q\in\LL^\times$. In view of the
factorization of $q$ given above, this is clearly equivalent with
the condition that $M=z^kM'$ holds for some $k\in\ZZ$. The projective
equivalence class of $M\in\Lat_n(\LL,\O_0)$ is denoted by
\[
[M]=\{z^kM|\ k\in\ZZ\}
\]
Thus we have obtained an action of the projective groups
$\PGL_n(\LL)$ and $\mathrm{PSL}_n(\LL)$ on the set 
\[
\{[M]|\ M\in\Lat_n(\LL,\O_0)\}
\]
of projective equivalence classes of $\O_0$-lattices.

\subsection*{The action of $\SL_n(\LL)$ and the type function}
The group $\SL_n(\LL)$ is \emph{not} transitive on set of projective
equivalence classes of $\O_0$-lattices. Let 
\[
M_0=\spa_{\O_0}\{e_1,\ldots,e_n\}
\]
denote the $\O_0$-module spanned by the canonical basis
$e_1,\ldots,e_n$ of $\LL^n$.
Suppose that $g(M_0)=M'$, for some $g\in\GL_n(\LL)$. Since
$\nu_0(\det(h))=0$ for all $h\in\GL_n(\O_0)$
(because $\nu_0(q)=0$ if $q\in\O_0$ is a unit),
the number
$\nu_0(\det(g))$ depends only on the module $M'$. The determinant of
the map $\lambda_{z^k}:v\mapstoo z^kv$ is $\det(\lambda_{z^k})=z^{kn}$.
Thus we have a well-defined map
\[
\Type([M'])=\nu(\det(g))+n\ZZ\in\ZZ/n
\]
which is $\SL_n(\LL)$-invariant. Note also that the stabilizers agree,
\[
\SL_n(\LL)_M=\SL_n(\LL)_{[M]},
\]
since $\det(\lambda_{z^k})\neq 1$ for $k\neq 0$. We put
\[
\V_i=\{[M]|\ M\in\Lat_n(\LL,\O_0), \Type([M])=i\}.
\]
Let
\[
M_i=\spa_{\O_0}\{ze_1,\ldots,ze_i,e_{i+1},\ldots,e_n\}
\]
denote the $\O_0$-module spanned by the vectors
$ze_1,\ldots,ze_i,e_{i+1},\ldots,e_n$. Then $[M_i]\in\V_i$.

\begin{Lem}
The action of $\SL_n(\LL)$ on $\V_i$ is transitive.

\proof
If $[M]=[g(M_i)]\in\V_i$ for some $g\in\GL_n(\LL)$, then
$\nu_0(\det(g))\equiv0\pmod n$, whence $\nu_0(\det(z^kg))=0$
for a suitable $k\in\ZZ$. Put $g'=z^kg$, then $[M]=[g'(M_i)]$,
and $\nu_0(\det(g'))=0$, whence $\det(g')\in\O_0$.
Finally, put $h=\mathrm{diag}(\det(g')^{-1},1,\ldots,1)$.
Then $h$ fixes $M_i$, and thus
$g'h\in\SL_n(\LL)$ maps $[M_i]$ to $[M]$.
\qed
\end{Lem}
These $n$ different $\O_0$-modules $M_0,\ldots,M_{n-1}$
thus form a cross-section for the
action of $\SL_n(\LL)$ on the set of projective equivalence classes of
$\O_0$-lattices, and we put
\begin{align*}
P^0_\rat=\SL_n(\LL)_{M_0}&=\SL_n(\O_0)\\
P^i_\rat=\SL_n(\LL)_{M_i}&=\left.\left\{\begin{pmatrix}A&z^{-1}B\\zC&D
\end{pmatrix}\right|\ A\in\O_0(i),\ 
\begin{pmatrix}A&B\\C&D\end{pmatrix}\in\SL_n(\O_0)\right\}
\end{align*}
for $i=1,\ldots,n-1$.
The sets $\V_0,\ldots,\V_{n-1}$ play the same r\^ole as the Grassmannians
$\Gr_j(\CC^n)$ in Section 1.1, and the stabilizers
$P_\rat^0,\ldots,P_\rat^{n-1}$ play the same r\^ole as the standard maximal
parabolics $P^1,\ldots,P^{n-1}$ in $\SL_n(\CC)$. However, there is one
fundamental difference:
the stabilizers $P_\rat^i$, $i=0,\ldots,n-1$, are conjugate to $P^0_\rat$
in $\GL_n(\LL)$,
\[
P_\rat^i=g_i\SL_n(\O_0) g_i^{-1},
\]
where $g_i=\mathrm{diag}(\underbrace{z,\ldots,z}_i,
\underbrace{1,\ldots,1}_{n-i})\in\GL_n(\LL)$.

\subsection*{Incidence, periodic flags, and apartments}
We define an \emph{incidence relation} $\tI$ on the set
$\{[M]|\ M\in\Lat_n(\LL,\O_0)\}$ as follows: 
\[
[M]\tI[M']\quad\overset{\text{def}}\Longleftrightarrow\quad
zM\leq z^kM'\leq M\text{ for some }k\in\ZZ .
\]
This relation is symmetric, since $zM\leq z^kM'\leq M$ implies that
$zM'\leq z^{1-k}M\leq M'$. Clearly, $\SL_n(\LL)$ acts by incidence
preserving automorphisms on this incidence geometry. Similarly as before,
we can use the incidence relation to construct a poset 
$\Delta(\LL^n,\O_0)$, the \emph{affine building} of $(\LL^n,\O_0)$;
the elements of $\Delta(\LL^n,\O_0)$ are the sets of pairwise
incident elements of $\V_0\cup\cdots\cup\V_{n-1}$.

Note also that $[M_i]\tI[M_j]$ holds for all $i,j$. The chambers are of the
following form. Let $\mathfrak B=(v_1,\ldots,v_n)$ be an ordered basis of
$\LL^n$, and put $M_i^{\mathfrak B}
=\spa_{\O_0}\{zv_1,\ldots,zv_i,v_{i+1},\ldots,v_n\}$.
Then $M({\mathfrak B})
=\{[M_0^{\mathfrak B}],\ldots,[M_{n-1}^{\mathfrak B}]\}$
is a maximal flag, which can be viewed as an infinite sequence of
free $\O_0$-modules
\[
\cdots
>z^{-1}M_{n-1}^{\mathfrak B}
>M_0^{\mathfrak B}
>M_1^{\mathfrak B}
>\cdots
>M_{n-1}^{\mathfrak B}
>zM_0^{\mathfrak B}
>zM_1^{\mathfrak B}
>\cdots
\]
of rank $n$. 
The quotient of two consecutive modules in this chain is a one-dimensional
complex vector space.
Note that $\Delta(\LL^n,\O_0)$ has rank $n$.

The collection of all chambers is called the \emph{periodic flag
variety} $\Fl(\LL^n,\O_0)$. The stabilizer of the chamber
$\{[M_0],\ldots,[M_{n-1}]\}$ is the Borel group
\[
B_\rat=P^0_\rat\cap P^1_\rat\cap\cdots\cap P^{n-1}_\rat=
\left\{
{\scriptscriptstyle\begin{pmatrix}\O_0&&\O_0\\&\ddots&\\z\O_0&&\O_0
\end{pmatrix}}\in\SL_n(\O_0)\right\}=ev_0^{-1}(\Sut_n(\CC)).
\]
%
%
Given a basis $v_1,\ldots,v_n$ of $\LL^n$, we
define the standard apartment $A^\LL_0\{v_1,\ldots,v_n\}$
as the collection of all partial flags obtained from the maximal flags
$M(\mathfrak B)$, where $\mathfrak B$ runs through the
collection of all bases of the form
$\mathfrak B=(z^{k_1}v_{\pi(1)},\ldots,z^{k_n}v_{\pi(n)})$,
for $\pi\in\Sym(n)$ and $k_1,\ldots,k_n\in\ZZ$.
The set-wise $\SL_n(\LL)$-stabilizer $N_\rat$ of
$A_0^\LL=A^\LL_0\{e_1,\ldots,e_n\}$
is the collection of all unimodular permutation matrices with entries in
$\LL^\times$, and the element-wise stabilizer $T_\rat$ of $A_0^\LL$ 
is the collection of all diagonal unimodular matrices with entries in
$\O_0^\times$. The quotient is the \emph{affine Weyl group}
$W=N_\rat/T_\rat\cong\tilde A_{n-1}$ of $\SL_n(\LL)$.
As a simplicial complex, $A_0^\LL$ is a triangulation of $\RR^{n-1}$.

\begin{Thm}
The simplicial complex $\Delta(\LL^n,\O_0)$ is a building of rank $n$
and type $\widetilde A_{n-1}$; as an apartment system, we may choose
the set $\cA^\LL_0=\{A^\LL_0\{v_1,\ldots,v_n\}|\ v_1,
\ldots,v_n\text{ a basis for }\LL^n\}$. The Coxeter diagram is
$\widetilde A_{n-1}$;
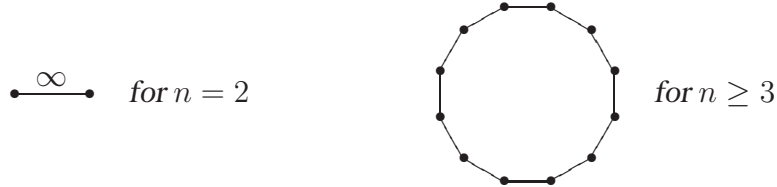
\[
\begin{xy} 
(0,0)*@{*};(10,0)*@{*} **@{-} ?(0.5)*!/_2mm/{\infty}
\end{xy}\quad\text{ for }n=2 \qquad\qquad\qquad
\begin{xy} /r12mm/:
 {\xypolygon12{@{*}}}
\end{xy}\quad\text{ for }n\geq 3 
\]
($n$ nodes).

\proof
Probably the easiest way to see that this is a building is to
verify that $(B_\rat,N_\rat)$ is a $BN$-pair for the group $\SL_n(\LL)$;
cp.~Brown \cite{Brown} Ch.~V.8 and Ronan \cite{Ronan} Ch.~9.2.
\qed
\end{Thm}
It is not difficult to see that the residue of $[M_i]$ is isomorphic to
$\Delta(\CC^n)$, for $i=0,\ldots,n-1$. Thus the panels in
$\Delta(\LL^n,\O_0)$ are again isomorphic to the complex projective
line $\CP^1$.

One final remark.
We have constructed the affine building $\Delta(\LL^n,\O_0)$
related to the discrete
valuation $\nu_0$. There is nothing special about the point $0\in\CP^1$;
if we choose a different point $c\in\CP^1$, then we obtain a different
affine building $\Delta(\LL^n,\O_c)$. These buildings are pairwise isomorphic;
in fact they are permuted by the group $\mathrm{PSL}_2(\CC)$.
Each of these building has a distinct collection of parabolics, so
$\SL_n(\LL)$ contains an uncountable set of $BN$-pairs.

\section{The twinning over $\CC[z,1/z]$}
\label{TwinBuildings}
Now we describe the buildings
$\Delta(\LL^n,\O_0)$ and $\Delta(\LL^n,\O_\infty)$
in a slightly different way, replacing the field $\CC(z)$ by the
ring $\CC[z,1/z]$. This section owes much to the paper \cite{AVM}
by Abramenko \& Van Maldeghem. So let
\[
\AA=\CC[z,1/z]=\bigcap\{\O_x|\ x\in\CC^\times\}
\]
denote the ring of all rational functions which are holomorphic
on $\CC^\times\subseteq\CP^1$. This is a subring of $\LL$, the ring of
\emph{Laurent polynomials}. Note that $\AA\cap\O_0=\CC[z]$ and
$\AA\cap\O_\infty=\CC[1/z]$. Similarly as before, we let
$\Lat_n(\AA,\CC[z])$ denote the collection of all free $\CC[z]$-modules
in $\AA^n$ which are spanned by an $\AA$-basis, and
$\Lat_n(\AA,\CC[1/z])$ the collection of all free $\CC[1/z]$-modules
spanned by $\AA$-bases.
Thus
\[
\Lat_n(\AA,\CC[z])=\{g(E_0^+)|\ g\in\GL_n(\AA)\} \quad\text{ and }\quad
\Lat_n(\AA,\CC[1/z])=\{g(E_0^-)|\ g\in\GL_n(\AA)\},
\]
where
$E_0^+=\spa_{\CC[z]}\{e_1,\ldots,e_n\}$ and 
$E_0^-=\spa_{\CC[1/z]}\{e_1,\ldots, e_n\}.$
For $E\in\Lat_n(\AA,\CC[z])\cup\Lat_n(\AA,\CC[1/z])$
we put as before
\[
[E]=\{z^kE|\ k\in\ZZ\},
\]
and
\[
E_i^+=\spa_{\CC[z]}\{ze_1,\ldots,ze_i,e_{i+1},\ldots,e_n\},
\quad
E_i^-=\spa_{\CC[1/z]}\{ze_1,\ldots,ze_i,e_{i+1},\ldots,e_n\},
\]
and
\[
\V_i^\pm=\{[gE_i^\pm]|\ g\in\SL_n(\AA)\} .
\]
The incidence is also defined as before.  If $[E]\in\V_i^\pm$ and
$[E']\in\V_j^\pm$, then
\[
[E]\tI[E']
\quad\overset{\text{def}}\Longleftrightarrow\quad
z^{\pm 1}E\leq z^k E'\leq E \quad\text{ for some }k\in\ZZ.
\]
Thus we obtain two simplicial complexes $\Delta^+(\AA^n)$ and
$\Delta^-(\AA^n)$
which are isomorphic under the map induced by the ring
automorphism $z\mapstoo 1/z$. The set of all maximal simplices
is denoted by $\Fl(\Delta^\pm(\AA^n))$.

Now there is a canonical map $\Lat_n(\AA,\CC[z])\too\Lat_n(\LL,\O_0)$
which maps $E$ to $\spa_{\O_0}(E)$, and a similar map
$\Lat_n(\AA,\CC[1/z])\too\Lat_n(\LL,\O_\infty)$; these maps induce
canonical $\SL_n(\AA)$-equivariant simplicial maps
$\Delta^+(\AA^n)\too\Delta(\LL^n,\O_0)$ and
$\Delta^-(\AA^n)\too\Delta(\LL^n,\O_\infty)$.
\begin{Prop}
The two maps $\Delta^+(\AA^n)\too\Delta(\LL^n,\O_0)$ and
$\Delta^-(\AA^n)\too\Delta(\LL^n,\O_\infty)$ are isomorphisms and thus
$\Delta^+(\AA^n)$, $\Delta^-(\AA^n)$ are buildings; the group
$\SL_n(\AA)$ acts transitively on both buildings.

\proof
For the proof we note that the $\SL_n(\AA)$-stabilizer of the 
$n-2$-simplex $B_i=\{[M_0],\ldots,\linebreak
{[M_{i-1}],[M_{i+1}],\ldots,[M_{n-1}]}\}
\in\Delta(\LL^n,\O_0)$ induces the transitive group $\mathrm{PSL}_2\CC$
on the corresponding panel. Thus $\SL_n(\AA)$ acts transitively on
$\Delta(\LL^n,\O_0)$ by Lemma \ref{TransLemma}.
Since $\SL_n(\AA)$ has the same stabilizers 
both in $\Delta^+(\AA^n)$ and in $\Delta(\LL^n,\O_0)$ (see below),
we obtain the claimed
isomorphism $\Delta^+(\AA^n)\cong\Delta(\LL^n,\O_0)$. The involution
$z\mapstoo 1/z$ on the ring $\AA$ and the projective line
$\CP^1=\CC\cup\{\infty\}$ normalizes $\SL_n(\AA)$, and we obtain
$\Delta^-(\AA^n)\cong\Delta(\LL^n,\O_\infty)$.
\qed
\end{Prop}
If $\mathfrak{B}=(v_1,\ldots,v_n)$ is an $\AA$-basis for $\AA^n$, then
we have similarly as before the chamber
$E^+(\mathfrak{B})=\{[E_0^{+,\mathfrak{B}}],\ldots,[E_{n-1}^{+,\mathfrak{B}}]\}$
corresponding to the modules
$E_i^{+,\mathfrak{B}}=\spa_{\CC[z]}\{zv_1,\ldots,zv_i,v_{i+1},\ldots,v_n\}$,
and the apartment $A^{+,\AA}\{v_1,\ldots,v_n\}$ obtained from the $\AA$-bases
$(z^{k_1}e_{\pi(1)},\ldots,z^{k_n}e_{\pi(n)})$; similarly, we obtain an
apartment $A^{-,\AA}\{v_1,\ldots,v_n\}$ for the building
$\Delta^-(\AA^n)$.
It can be checked that $\SL_n(\AA)$ acts strongly transitively on
$\Delta^+(\AA^n)$ and $\Delta^-(\AA^n)$ with respect to these apartment
systems
\[
\cA^\pm=\{A^{\pm,\AA}\{v_1,\ldots,v_n\}|\ v_1,\ldots,v_n
\text{ an $\AA$-basis for }\AA^n\}.
\]
Note that the apartment system $\cA^{+,\AA}$ is strictly smaller than
$\cA^\LL_0$.
\begin{Prop}
The group $\SL_n(\AA)$ acts strongly transitively both on $\Delta^+(\AA^n)$
and on $\Delta^-(\AA^n)$.

\proof This follows from the fact that $\SL_n(\AA)$ acts transitively on
pairs of \emph{opposite} chambers in $\Delta^+(\AA^n)\times\Delta^-(\AA^n)$,
cp.~Abramenko \& Van Maldeghem \cite{AVM};
the opposition relation is defined in the next section.
\qed
\end{Prop}
We put
\[
P_\alg^{i,+}=P_\rat^i\cap\SL_n(\AA)
\quad\text{ and }\quad
B_\alg^+=B_\rat\cap\SL_n(\AA).
\]
Thus
\[
P^{0,+}_\alg=\SL_n(\CC[z])\quad\text{ and }\quad
P^{i,+}_\alg=g_i(P^{0,+}_\alg)g_i^{-1},
\]
where $g_i=\mathrm{diag}(z,\ldots,z,1,\ldots,1)$ as before.
These are the parabolics corresponding to $\Delta^+(\AA^n)$; 
there are similar parabolics for $\Delta^-(\AA^n)$.

\subsection*{Twin buildings}

The notion of a twin building was developed by Ronan and Tits in order
to supply geometries for Kac-Moody groups; we refer to
Abramenko \cite{AbraLNM}, Abramenko \& Ronan \cite{AR}
Abramenko \& Van Maldeghem~\cite{AVM},
M\"uhlherr \cite{Mu3} \cite{Mu2}, 
M\"uhlherr \& Ronan \cite{MR},
Ronan \cite{Ron2}, Ronan \& Tits \cite{RT}
and Tits \cite{Durham},
\cite{TitsRes} for more details about twin buildings.

Let $(W,S)$ be a Coxeter system, and let $(\Delta^+,\Delta^-)$ be a
pair of buildings with this given Coxeter system. The $W$-valued distance
in $\Delta^\pm$ is denoted $\delta_{\Delta^\pm}$.
A \emph{twinning} of $\Delta^+$ with $\Delta^-$ is 
a $W$-valued \emph{codistance} function
\[
\delta^*:\Cham(\Delta^\pm)\times\Cham(\Delta^\mp)\too W,
\]
subject to the following axioms. 
\emph{The intuitive idea is that objects
with a small codistance are far apart.}

\smallskip\textbf{Tw}$_1$ 
The relation
$\delta^*(C^\pm,D^\mp)=\delta^*(D^\mp,C^\pm)^{-1}$ holds for all
chambers $C^\pm\in\Delta^\pm$, $D^\mp\in\Delta^\mp$
(the codistance is 'symmetric').

\smallskip\textbf{Tw}$_2$
Let $w\in W$ and $s\in S$, and suppose that $\ell(ws)=\ell(w)-1$
(i.e. that $w$ has a reduced expression with $s$ as the last letter).
If $C^\pm\in\Delta^\pm$, $D^\mp,E^\mp\in\Delta^\mp$ are chambers
with $\delta^*(C^\pm,D^\mp)=w$ and $\delta_{\Delta^\mp}(D^\mp,E^\mp)=s$,
then $\delta^*(C^\pm,E^\mp)=ws$
(all chambers $E^\mp$ in the $s$-panel through $D^\mp$ are
'further away' from $C^\pm$).

\smallskip\textbf{Tw}$_3$
If $C^\pm\in\Delta^\pm$, $D^\mp\in\Delta^\mp$ are chambers with codistance
$\delta^*(C^\pm,D^\mp)=w$, and if $s\in S$, then there exists a
chamber $E^\mp\in\Delta^\mp$, with $\delta_{\Delta^\mp}(D^\mp,E^\mp)=s$,
and with $\delta^*(C^\pm,E^\mp)=ws$ (the codistance leads to galleries).

\smallskip
Note that there is a symmetry in the axioms if we exchange the signs
'$+$' and '$-$'. Sometimes we will state a result for a specific choice
of the signs; the corresponding result for the opposite choice of signs
follows in the same way.

It is clear how to extend the codistance to a double coset valued
codistance
\[
\delta^*:\Delta^\pm\times\Delta^\mp\too
\textstyle\bigcup\{W_J\backslash W/W_K|\ J,K\subseteq I\}.
\]
The \emph{opposition relation}
$\op\subseteq \Delta^+\times\Delta^-$ is defined as follows:
two chambers are opposite if their codistance is 1; this relation
extends in a natural way to the simplices. The codistance can be used
to sync the type functions in both buildings in such a way that
opposite vertices have the same type, and we will assume
this to be done.

A (special) \emph{automorphism} of a twin building is a pair $(g^+,g^-)$
of automorphisms $g^\pm\in\Spe(\Delta^\pm)$ which preserves the codistance,
$\delta^*(C^+,C^-)=\delta^*(g^+(C^+),g^-(C^-))$.

\begin{Num}
\textbf{Example}
Let $\Delta$ be a spherical building, let $w_0$ denote the unique longest
element in the Coxeter group $W$ of $\Delta$, put $\Delta^+=\Delta^-=\Delta$,
and $\delta^*(C,D)=\delta(C,D)w_0$. The resulting geometry is a twin
building.
\end{Num}
All twin buildings with spherical halves arise in this way, see
Tits \cite{TitsRes}. Twin buildings
are natural generalizations of spherical buildings; they share
many of the particular geometric properties of spherical buildings.
The twinning is in general not determined by the pair $(\Delta^+,\Delta^-)$;
a pair of buildings (e.g.~a pair of trees) can admit many
non-equivalent twinnings.

The group $\SL_n(\AA)$ induces a twinning on the pair
$(\Delta^+(\AA^n),\Delta^-(\AA^n))$;
the codistance can be defined as follows. The group $\SL_n(\AA)$
has a \emph{Birkhoff decomposition} (a Bruhat twin decomposition) as
\begin{align*}
\SL_n(\AA)&=B^-_\alg N B^+_\alg \\
&=\textstyle\dot{\bigcup}\{B^-_\alg w B^+_\alg|\ w\in N/
(B^+_\alg\cap B^-_\alg)\} 
\end{align*}
where $N$ is the set-wise
stabilizer of the twin apartment
$(A^{+,\AA}\{e_1,\ldots,e_n\},A^{-,\AA}\{e_1,\ldots,e_n\})$. Since
$\Cham(\Delta^\pm(\AA^n))=\SL_n(\AA)/B^\pm_\alg$ and
$N/(B^+_\alg\cap B^-_\alg)\cong W$, we may use the
Birkhoff decomposition to define the codistance as
\[
\delta^*(g B^-_\alg,h B^+_\alg)=w
\quad\text{ if and only if }\quad
B^-_\alg g^{-1}hB^+_\alg=B^-_\alg w B^+_\alg.
\]
The following Theorem is 'folklore'. A nice proof is given in
Abramenko \& Van Maldeghem \cite{AVM}.
\begin{Thm}
The triple $(\Delta^+(\AA^n),\Delta^-(\AA^n),\delta^*)$
is a twin building, and the
group $\SL_n(\AA)$ acts as a strongly transitive group of automorphisms.
Two chambers $(C^+,C^-)\in\Delta^+(\AA^n)\times\Delta^-(\AA^n)$ are opposite if and
only if they arise from lattices which 
are 'back to back', i.e. if there exists an ordered $\AA$-basis
$(v_1,\ldots,v_n)$ such that $C^+=E^+(v_1,\ldots,v_n)$ and
$C^-=E^-(v_1,\ldots,v_n)$.
\end{Thm}

\section{Loop groups}

There is one big difference between $\SL_n(\LL)$ and $\SL_n(\AA)$
which is due to the fact that $\AA$ is a ring, while $\LL$ is a field:
the group $\SL_n(\LL)$ is almost simple, whereas the
group $\SL_n(\AA)$ is a semi-direct product. 

\subsection*{Based loops}

Fix $c\in\CC^\times$ and
consider the evaluation map
$ev_c:\SL_n(\AA)\too\SL_n(\CC)$ given by
\[
ev_c:f(z)\mapstoo f(c).
\]
We denote the kernel of this map by $\Omalg(\SL_n(\CC),c)$
and we put $\Omalg\SL_n(\CC)=\Omalg(\SL_n(\CC),1)$ for short.
The injection $\SL_n(\CC)\subseteq\SL_n(\AA)$ leads to a split exact
sequence
\[
1\too\Omalg(\SL_n(\CC),c)\too\SL_n(\AA)\too^{ev_c}_{\leftarrow}\SL_n(\CC)
\too1,
\]
hence $\SL_n(\AA)$ is a semi-direct product.
(There is an obvious $\CC^\times$-action on $\SL_n(\AA)$
given by $f(z)\mapstoo f(az)$; under this action, the collection
of normal subgroups $\{\Omalg(\SL_n(\CC),c)|\ c\in\CC^\times\}$ is 
permuted transitively.
Thus one is lead to the semi-direct product $\SL_n(\AA)\rtimes\CC^\times$.)
\begin{Lem}
For every $c\in\CC^\times$, the group $\Omalg(\SL_n(\CC),c)$ acts
transitively on each of the sets
$\V_0^+,\V_0^-,\ldots,\V_{n-1}^+,\V_{n-1}^-$.

\proof
We clearly have $\SL_n(\CC)\subseteq\SL_n(\CC[z])=P_\alg^{0,+}$, and thus
\[
P_\alg^{0,+}\,\Omalg(\SL_n(\CC),c)\supseteq\SL_n(\CC)\,
\Omalg(\SL_n(\CC),c)=\SL_n(\AA).
\]
Now $P_\alg^{i,+}=g_iP_\alg^{0,+}g_i^{-1}$ and
$\Omalg(\SL_n(\CC),c)$ is $g_i$-invariant, whence
\[
P_\alg^{i,+}\Omalg(\SL_n(\CC),c)=g_iP_\alg^{0,+}g_i^{-1}\Omalg(\SL_n(\CC),c)=
\SL_n(\AA).
\]
\qed
\end{Lem}
There is a natural map
$\SL_n(\AA)\too C^\infty(\SS^1,\GL_n(\CC))=L_\diff\SL_n(\CC)$
into the set $L_\diff\SL_n(\CC)$
of smooth maps from the unit circle into $\SL_n(\CC)$; this map
is obtained by viewing the elements of $\SL_n(\AA)$
as maps from $\SS^1\subseteq\CC^\times$ into $\SL_n(\CC)$. If
$f=\sum_\fin f_kz^k\in\SL_n(\AA)$, then 
\[
f_k=\frac{1}{2\pi\tti}\oint_{|z|=1}\frac1{z^{k+1}}f(z) dz.
\]
Therefore, the map into $L_\diff\SL_n(\CC)$ is an injection. From now
on, we denote the group $\SL_n(\AA)$ also by $L_\alg\SL_n(\CC)$.
The subgroup $\Omalg(\SL_n(\CC),c)$ is thus the subgroup
$L_\alg\SL_n(\CC)\cap\Omega_\diff(\SL_n(\CC),c)$ of $c$-based 
algebraic loops
(by $c$-based we mean that the base-point of $\SS^1$ is
$c$ --- the base point of the group is always the identity element).

\subsection*{The Cartan involution}
Consider the semi-linear involution $\#$ on the complex
vector space $\AA(n)$ which is given by
\[
f=\sum_\fin f_kz^k\mapstoo f^\#=\sum_\fin f^*_{-k}z^k.
\]
The subgroup of all elements $g\in\SL_n(\AA)$ with $g^{-\#}=g$
is denoted by $L_\alg\SU(n)$. This terminology is motivated by the
fact that
\[
L_\alg\SU(n)=L_\diff\SU(n)\cap L_\alg\SL_n(\CC).
\]
To see this, note that for $z\in\SS^1$ we have $z^{-1}=\bar z$, hence
$f(z)f(z)^*=1$ holds for all $z\in\SS^1$ if and only if
$f\in L_\alg\SU(n)$. For $c\in\SS^1$ we have a semi-direct
decomposition
\[
L_\alg\SU(n)=\SU(n)\Omalg(\SU(n),c)
\]
as before, and an $\SS^1$-action which is given by
$f(z)\mapstoo f(az)$. 

\begin{Lem}
We have $L_\alg\SU(n)\cap\SL_n(\CC[z])=\SU(n)$.

\proof
Let $f=f_0+f_1z+\cdots+f_kz^k\in\SL_n(\CC[z])$. Then
$f^\#=f^*_kz^{-k}+\cdots f^*_1z^{-1}+f_0^*$. If
$f^\#=f^{-1}\in\SL_n(\CC[z])$, then $f_i=0$ for $i\geq 1$.
\qed
\end{Lem}
In particular,
\[
B_\alg^+\cap L_\alg\SU(n)\cong T^{n-1}
\quad\text{ and }\quad
P^{i,+}_\alg\cap L_\alg\SU(n)\cong\mathrm{S}(\U(i)\cdot\U(n-i))
\]
(as in Section \ref{SphericalCase},
$T^k$ denotes a compact torus of rank $k$).
There is a more geometric description of the group $L_\alg\SU(n)$.
The map $\AA^n\too^{\#}\AA^n$ induces isomorphisms
$\Delta^+(\AA^n)\too^{\#}\Delta^-(\AA^n)$ and
$\Delta^-(\AA^n)\too^{\#}\Delta^+(\AA^n)$,
\[
\begin{xy}
{(0,0)*{\Delta^+(\AA^n)\qquad}="a"}
,{(20,0)*{\qquad\Delta^-(\AA^n)}="b"}
,{\ar@{<->}@/^3ex/"a";"b"}
,{"a";"b"**@{} ?(0.5)*!/_9mm/{\scriptstyle\#}}
\end{xy}
\]
and $L_\alg\SU(n)=\mathrm{Cen}_{\SL_n(\AA)}(\#)$.
There is a corresponding Cartan decomposition of the loop algebra
$\fsl_n(\AA)$ into eigenspaces of the involution $\#$,
\[
\fsl_n(\AA)=L_\alg\mathfrak{su}(n)\oplus\mathfrak{X},
\]
where $\mathfrak X$ denotes the traceless hermitian matrices in $\AA(n)$.
\begin{Thm}
The group $L_\alg\SU(n)$ acts transitively on the periodic flags,
\[
L_\alg\SL_n(\CC)=B_\alg^+ L_\alg\SU(n),
\qquad
\Fl(\Delta^+(\AA^n))\cong L_\alg\SU(n)/T^{n-1}.
\]
\proof
The proof is exactly the same as in \ref{SU(n)IsTransitive}.
The $L_\alg\SU(n)$-stabilizer of a panel
is isomorphic to $\SU(2)\cdot T^{n-2}$ and induces the transitive group
$\SO(3)$ on the panel.
\qed
\end{Thm}
In particular, $L_\alg\SU(n)$ acts transitively on $\V_0^+$. Since
$L_\alg\SU(n)_{[M_0]}=\SU(n)$, the group $\Omalg\SU(n)$ acts
\emph{regularly} on $\V_0^+$,
\[
\V_0^+\cong L_\alg\SU(n)/\SU(n)\cong\Omalg\SU(n),
\]
and
\[
\Fl(\Delta^+(\AA^n))\cong\Omalg\SU(n)\times(\SU(n)/T^{n-1}).
\]
Note that the proof above
(which is due to Bernhard M\"uhlherr) is much simpler
than the classical proof given e.g.~in Pressley \& Segal 
\cite{PresSeg} Theorem 8.3.2.

\section{The affine Veronese representation of $\Delta^+(\AA^n)$}

In this section we construct an equivariant embedding
$\Delta^+(\AA^n)\rInto \fsl_n(\AA)$ which is very similar to the
finite-dimensional Veronese representation
$\Delta(\CC^n)\rInto \fp_n\subseteq\fsl_n(\CC)$.
Recall that we presented the flags in $\Delta(\CC^n)$ as certain
hermitian operators. Similarly, we want to associate an operator to
the $\CC[z]$-module $E_0$. To this end we
consider the first order linear differential operator
\[
\AA\too\AA,
\qquad f\mapstoo \zpz f=z\frac{\partial f}{\partial z}.
\]
If $f=\sum_\fin f_kz^k$, then $\zpz f=\sum_\fin kf_kz^k$.
Thus 
\[
\kernel(\zpz-\lambda)=
\begin{cases}
z^\lambda\CC &\text{ for }\lambda\in\ZZ \\
0& \text{ for }\lambda\in\CC\setminus\ZZ .
\end{cases}
\]
In particular,
\[
\CC[z]=\bigoplus_{k\geq 0}\kernel(\zpz-k)
\quad\text{ and }\quad
\CC[1/z]=\bigoplus_{k\geq 0}\kernel(\zpz+k).
\]
The operator $\zpz$ extends in an obvious way to $\AA^n$ and
to the matrix algebra $\AA(n)$. For $f\in\AA^n$ we put
\[
Df=\zpz f .
\]
Let $\Pi_k=\mathrm{diag}(\underbrace{1,\ldots,1}_k,0,\ldots,0)\in H(n)$
and let $\Pi_k^\tl=\Pi_k-\frac kn\one$ denote its traceless image in
$\fp_n$.
Thus
\begin{align*}
E_i^+&=\spa_{\CC[z]}\{ze_1,\ldots,ze_i,e_{i+1},\ldots,e_n\} \\
&=\bigoplus_{k\geq0}\kernel(D-\Pi_i-k\one) \\
&=\bigoplus_{k\geq0}\textstyle\kernel((D-\Pi_i^\tl)-(k+\frac in)\one);
\end{align*}
the elements of the flag varieties $\V_i^+$ correspond
bijectively to the
$L_\alg\SU(n)$-conjugates of $D-\Pi_i^\tl$. Let $g\in L_\alg\SU(n)$.
Then
\[
0=\zpz(gg^\#)=(\zpz g)g^\#+g\zpz(g^\#).
\]
Therefore
$g D g^\#f=g(\zpz(g^\#))f+gg^\#Df=Df-(\zpz g)g^\# f$, whence
$gDg^\#=D-(\zpz g)g^\#$, and
\[
g(D-\Pi_i^\tl)g^\#=D-(\zpz g)g^\#-g\Pi_i^\tl g^\#.
\]
We put 
\[
\fX=\{X\in\fsl_n(\AA)|\ X^\#=X\}.
\]
This is an infinite-dimensional
real vector space, and $(z\partial_zg)g^\#\in\fX$ for all $g\in L_\alg\SU(n)$.
We construct an
$L_\alg\SU(n)$-equivariant embedding of $\Delta^+(\AA^n)$ into the
infinite-dimensional real vector space $\fX\oplus \RR D$ as follows.
\begin{align*}
\V_k^+ & \rInto \fX\oplus \RR D \\
[gE_k^+] & \mapstoo g(D-\Pi_k^\tl)g^\# = D-(\zpz g)g^\#-
g\Pi_k^\tl g^\#.
\end{align*}
We extend this mapping to the geometric realization $|\Delta^+(\AA^n)|$
of $\Delta^+(\AA^n)$ in the canonical way, and we call the resulting
$L_\alg\SU(n)$-equivariant map 
\[
|\Delta^+(\AA^n)|\rInto \fX\oplus\RR D\too^{-pr_1}\fX
\]
the \emph{affine Veronese representation} of $\Delta^+(\AA^n)$.
Explicitly, the Veronese representation of the vertex $[gE_k]$ is
\[
\Phi([gE_k])=g\Pi_k^\tl g^\#+(\zpz g)g^\#.
\]
Note that the group $L_\alg\SU(n)$ acts through \emph{gauge transformations}
\[
X\mapstoo gXg^\#+(\zpz g)g^\#
\]
on $\fX$.
Similarly as in the spherical case, it is not difficult to prove that
$\Phi$ injects the geometric realization $|\Delta^+(\AA^n)|$ into
$\fX$.
The partial flags in the building $\Delta^+(\AA^n)$ correspond
thus to certain operators $D-X\in \fX\oplus\RR D$ with finite-dimensional
nontrivial eigenspaces. Note also that $\#$ swaps $\Delta^+(\AA^n)$ and
$\Delta^-(\AA^n)$, so $\Phi$ is at the same time a Veronese
representation for $\Delta^-(\AA^n)$ with exactly the same image.

\begin{Num}\textbf{Caveat}
The affine Veronese representation
is \emph{not} surjective. Let $a(z)=z+1/z$ and $r\in\RR^\times$.
The differential equation
\[
(\zpz -ra)f=\lambda f
\]
has the solution $f(z)=e^{r(z-1/z)}z^\lambda\cdot\mathrm{const}$.
This function is holomorphic on $\CC^\times$ if and only if $\lambda$
is an integer, but it is not meromorphic on $\CC\cup\{\infty\}$,
so $f\not\in\AA$. Thus, the traceless
diagonal matrix $X=\mathrm{diag}(a,\cdots,a,(1-n)a)\in\fX$
does \emph{not} represent a flag of the building
because $D-X$ has no nontrivial eigenspaces.
\end{Num}
The space $\fX$ is a subspace of the \emph{loop algebra}
$\fsl_n(\AA)$, which in turn is contained in the semi-direct
product $\fsl_n(\AA)\oplus\CC D$. So $\fX\oplus \RR D$ plays a very similar
r\^ole as $\fp_n\subseteq\fsl_n(\CC)$. Let $Q$ denote the collection
of all barycenters of images of chambers in $\fX$. The closure $M$ of $Q$
in the Hilbert space completion of the real pre-Hilbert space $\fX$ is
an \emph{infinite-dimensional isoparametric submanifold}.
See Pinkall \& Thorbergsson \cite{PT} and
Heintze, Palais, Terng \& Thorbergsson \cite{HPTT} for more
information. 
The set $Q$ coincides with the subset $Q(p)\subseteq M$ introduced in 
\cite{HeintzeLiu} by Heintze \& Liu.
In Gro\ss, Heintze, Kramer \& M\"uhlherr \cite{GHKM}
we prove that all known isoparametric submanifolds of rank
at least 3 in Hilbert spaces
arise in a uniform way from Veronese representations of
twin buildings.

\section{Topological geometry and Bott periodicity}
\label{TopologicalGeometry}
In this section we propose a definition of topological twin buildings.
Since spherical buildings are twin buildings, this is at the same time
a definition of topological spherical buildings. 
Definitions of spherical topological buildings have been proposed by
Burns \& Spatzier \cite{BurnsSpatzier}, J\"ager \cite{Jae},
K\"uhne \cite{Kue}, and myself; Mitchell \cite{Mitchell} proposes
an \emph{ad hoc} definition of topological $BN$-pairs.
For spherical buildings of rank 2 there is a well-established theory, see
Salzmann \cite{Salz}, Salzmann \emph{et al.} \cite{CPP},
Schroth \cite{Schroth},
Grundh\"ofer \& Knarr \cite{GK}, 
Grundh\"ofer \& Van Maldeghem \cite{GVM},
Grundh\"ofer, Knarr \& Kramer \cite{Flaghom1}, 
and Kramer \cite{KraDiss}, \cite{KraHabil}.
The starting point is always
a topology on the set of vertices of the building
(i.e. on the 0-simplices). Using the type function,
the simplices of higher rank can be
interpreted as ordered tuples of vertices, and thus one obtains a topology
on $\Delta$; the question then is which maps should be
continuous. Burns \& Spatzier \cite{BurnsSpatzier} require only that
$\Cham(\Delta)$ should be closed, i.e. that every net of chambers,
viewed as a net of $r$-tuples of vertices, converges to some chamber.
Moreover, they claim that for the classical geometries (projective
spaces) this agrees with the traditional notion of
a topological geometry, see \emph{loc.cit}~p.~1.
This is definitely not true, and it is not difficult to construct
perverse topologies on nice
geometries which satisfy their condition nevertheless.
However, in the \emph{compact} spherical case, their definition is the
correct one (and that's the only instance where they need it in their
work \cite{BurnsSpatzier}), see Proposition \ref{TheoProp} below.

Compactness or local compactness in the non-spherical case leads to
locally finite buildings (finite panels), and this in turn is related
to locally compact CAT(0)-spaces; however, these matter are not within
the scope of the present article.

Similarly as in J\"ager \cite{Jae}, our notion of a topological building
asks for the continuity of certain \emph{projections}. The fact that
this definition makes sense for twin buildings was pointed
out by Bernhard M\"uhlherr during a meeting in Oberwolfach 
back in 1992; then, we planned to write a joint paper with
Martina J\"ager on projections and
topologies in buildings which, however, never came to existence.
This section is a first approximation of what we
had in mind.

The fact that topological buildings can be used to prove Bott periodicity
is particularly appealing, since topological $K$-theory is an important
ingredient in topological geometry; many classification results in 
Salzmann \emph{ et al.} \cite{CPP} and in Kramer \cite{KraDiss}
depend in an essential way on Bott periodicity.

\subsection*{Projections in (twin) buildings}
Tits defined projections in (spherical) buildings in \cite{Tits};
a modern account based on metric properties of buildings is
given in Dress \& Scharlau \cite{DS}.
Let $C,D$ be chambers in a building $\Delta$, let $\delta(C,D)=w$, and
let $w=s_{i_1}\cdots s_{i_r}$ be a reduced (minimal) expression for
$w$ in terms of the generating set $S$. Then there exists a unique
\emph{minimal gallery} $\gamma=(C_0,C_1,C_2,\ldots,C_r=D)$ 
\emph{of type} $(s_{i_1},\ldots,s_{i_r})$, consisting of chambers
$C_0,\ldots,C_r$, such that $\delta(C_{k-1},C_k)=s_{i_k}$ holds
for $k=1,\ldots,r$.
Now let $X\in \Delta$ be a simplex and let $C$ be a chamber.
Then there exists a unique chamber $E$ in $\Res(X)$ which we denote
\[
E=\proj_XC,
\]
the \emph{projection of $C$ onto $X$},
with the following
property: for every chamber $D\in\Res(X)$, and for every minimal
gallery $\gamma$ starting at $C$ and ending at $D$, the first chamber
in $\gamma$ which is contained in $\Res(X)$ is $E$, the \emph{gate}
of $\Res(X)$ with respect to $C$.
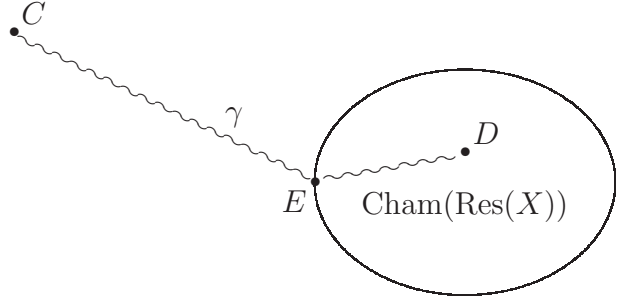
\[
\begin{xy}/r2cm/:0,
{(0,1)*\dir{*}="c"}
,{*+!DL{C}}
,{(2,0)*\dir{*}="e"}
,{*+!UR{E}}
,{(3,0.2)*\dir{*}="d"}
,{*+!DL{D}}
,{(3,0)="m"}
,{*+!Ul{\Cham(\Res(X))}}
,"e","m",{\ellipse(1,.75){}}
,"c";"e"**@{~} ?(0.7)*!/_3mm/{\gamma}
,"e";"d"**@{~}
\end{xy}
\]
Note that $\proj_XC=C$ if $C\in\Cham(\Res(X))$. If $Y\in\Delta$ is an
arbitrary simplex, then there exists a unique simplex
$Z$ which is contained in some chamber in $\Res(X)$,
such that
\[
\Cham(\Res(Z))=\proj_X\Cham(\Res(Y)),
\]
and we put
\[
Z=\proj_XY.
\]
Now suppose that $(\Delta^+,\Delta^-,\delta^*)$ is a twin
building, that $X\in\Delta^+$ is a \emph{spherical} simplex
(recall from Section 1 that this means that $\Res(X)$ is spherical), and 
that $C\in\Delta^-$ is a chamber. Then there exists
a unique chamber $E\in\Res(X)$ which maximizes the numerical
codistance function $D\mapstoo\ell(\delta^*(C,D))$ on $\Cham(\Res(X))$,
see Ronan \cite{Ron2} (4.1).
Intuitively, a 'small' codistance corresponds to a 'big'
distance, so $E$ is the chamber 'closest' to $C$; note that
$D\mapstoo\ell(\delta^*(C,D))$ is bounded above because $X$ is spherical.
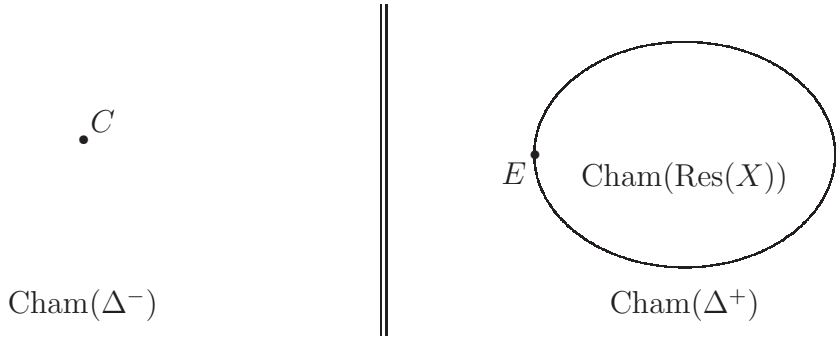
\[
\begin{xy}/r2cm/:0,
{(-1,.1)*\dir{*}="c"}
,{*+!DL{C}}
,{(2,0)*\dir{*}="e"}
,{*+!UR{E}}
,{(-1,-1)*{\Cham(\Delta^-)}}
,{(3,-1)*{\Cham(\Delta^+)}}
,{(3,0)="m"}
,{*+!Ul{\Cham(\Res(X))}}
,"e","m",{\ellipse(1,.75){}}
,(1,-1.2);(1,1)**@{=}
\end{xy}
\]
The chamber $E$ is again denoted
\[
E=\proj_XC,
\]
and the projection $\proj_XY$ of an arbitrary simplex $Y\in\Delta^-$ onto
$X$ is defined exactly as before (but only for spherical $X$!).

\begin{Num}\
\textbf{Schubert cells}
Let $C_0$ be a chamber in a building $\Delta$, let $J\subseteq I$,
and let $w W_J\in W/W_J$. The set 
\[
\cC_{wW_J}(C_0)=\{X\in\Delta|\ \delta(C_0,X)=wW_J\}
\]
is called a \emph{Schubert cell} in $\Delta$. Schubert cells
in general buildings have no special structure (e.g.~if $\Delta$ is a tree),
but Schubert cells in halves of twin buildings (and in particular
Schubert cells in spherical buildings) have a nice product structure,
i.e. they admit coordinates (labels).
\end{Num}
\begin{Num}
\textbf{Coordinatizing a half twin}
\label{Coordinatization}
Let $C_0\in\Delta^+$ and $D_0\in\Delta^-$ be a pair of opposite chambers
in a twin building $(\Delta^+,\Delta^-,\delta^*)$.
These two chambers determine a unique apartment $A\subseteq\Delta^-$
(half of the \emph{twin apartment} spanned by $(C_0,D_0)$, see
e.g.~Ronan \cite{Ron2} 2.8).
Let $w\in W$, and let $w=s_{i_1}\cdots s_{i_r}$ be a reduced expression.
Let $E\in \cC_w(C_0)$, and let $C_0,\ldots,C_r=E$ be a 
(necessarily minimal) gallery of type $(s_{i_1},\ldots,s_{i_r})$.
Let $D_0,\ldots,D_r$ be the unique minimal gallery of the same type
in the apartment $A\subseteq\Delta^-$.
Then $C_k\,\op\,D_k$ holds for 
$k=0,\ldots,r$. We define $r$ \emph{coordinates}
$(X_1,\ldots,X_r)$ by $X_k=\proj_{D_k\sqcap D_{k-1}}C_k$.
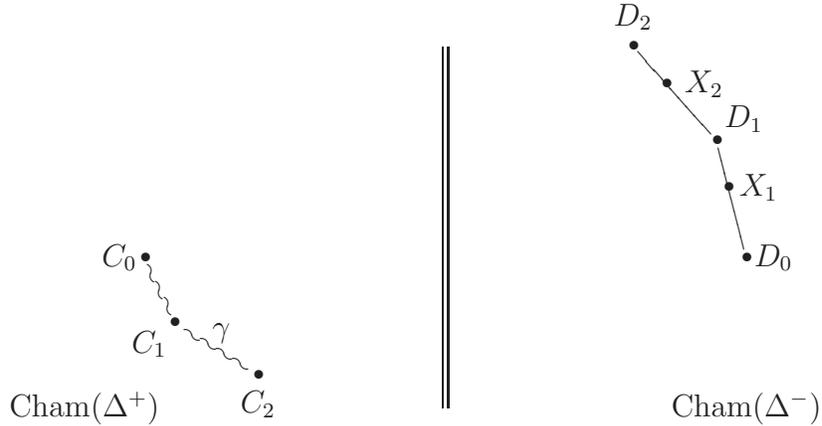
\[
\begin{xy}/r4cm/:
0
,{\xypolygon14"A"{~:{(1,0):(0,0.9)::}~*{}~>{}}}
,{\xypolygon14"B"{~:{(1,0):(0,0.5)::}~*{}~>{}}}
,(0,-0.5);(0,0.7)**@{=}
,"A1"*+!L{D_0}
,"A2"*+!LD{D_1}
,"A3"*++!D{D_2}
,"B8"*+!R{C_0}
,"B9"*+!RU{C_1}
,"B10"*++!U{C_2}
,"A1"*@{*}
,"A2"*@{*}
,"A3"*@{*}
,"A1";"A2"**@{-}?(0.6)*@{*}
,"A1";"A2"**@{}?(0.6)*!/l4mm/{X_1}
,"A2";"A3"**@{-}?(0.6)*@{*}
,"A2";"A3"**@{}?(0.6)*!/l5mm/{X_2}
,"B8"*@{*}
,"B9"*@{*}
,"B10"*@{*}
,"B8";"B9"**@{~}
,"B9";"B10"**@{~}?(0.2)*!/l4mm/{\gamma}
,{(-1.2,-0.5)*{\Cham(\Delta^+)}}
,{(1,-0.5)*{\Cham(\Delta^-)}}
\end{xy}
\]
Note that $X_k\in\Res(D_k\sqcap D_{k-1})\setminus\{D_k\}$.
Now the point is that step by step, the whole gallery $\gamma$ can be
recovered from these coordinates;
\begin{align*}
C_1&=\proj_{C_0\sqcap C_1}X_1 \\
C_2&=\proj_{C_1\sqcap C_2}X_2 \\
C_3&=\proj_{C_2\sqcap C_3}X_2 \quad \text{\em etc.}
\end{align*}
since $C_{k-1}$ and $s_{i_k}$ determine $C_{k-1}\sqcap C_k$,
the information needed is only $C_0$, the coordinates $(X_1,\ldots,
\newline X_r)$,
and the reduced expression $(s_{i_1},\ldots,s_{i_r})$. Note also that
different reduced expressions for a Schubert cell lead to different
coordinates; our coordinatization process depends on a choice of a
reduced expression for every element $w\in W$. In the case of
spherical buildings of rank two, the different expressions of
the longest element in the Coxeter group and the resulting
'changes of coordinates' lead to Van Maldeghem's coordinatizing
rings \cite{HVM}.
\end{Num}
The process above yields coordinates for the Schubert cells
$\cC_w(C_0)\subseteq\Cham(\Delta^+)$.
For a coset $wW_J\in W/W_J$ we may assume that $w$ is the unique
\emph{shortest coset representative}. The canonical 'forgetful' map
$\cC_w(C_0)\too\cC_{wW_J}(C_0)$ which sends a chamber $C$ to the unique
subsimplex $C|_{I\setminus J}$ of type $I\setminus J$ contained in $C$
is a bijection; in fact, $C=\proj_{C|_{I\setminus J}}C_0$.
The Schubert cell $\cC_{wW_J}(C_0)$
can also directly be coordinatized, by the same
method as described above.
In any case we see that each Schubert cell is in a natural correspondence
with a finite product of punctured panels; for $\Delta^+(\AA^n)$ we
see that the Schubert cell $\cC_{wW_J}(C_0)$ bijects onto $\CC^m$, 
where $m$ is the length of the shortest coset
representative of $wW_J$.

\subsection*{Topological twin buildings}

Let $(\Delta^+,\Delta^-,\delta^*)$ be a twin building. Suppose that
there is a Hausdorff topology on the set of vertices
(the 0-simplices) of both buildings. The simplices of type $J$ can be
regarded as $J$-tuples of vertices; in this way, the topology on the
vertices determines a topology on both buildings.
For $J,K\subseteq I$ and $w\in W$ we put
\[
\cD_{W_JwW_K}^{J,K}=
\{(X,Y)\in\Delta^+\times\Delta^-|\ 
\Type(X)=I\setminus J,\,
\Type(Y)=I\setminus K,\,
\delta^*(X,Y)=W_JwW_K\}.
\]

\begin{Def}
A twin building
is called a \emph{topological twin building} if the following condition
is satisfied:

\smallskip\textbf{TTB}
If $J\subseteq I$ is spherical (i.e. if $W_J$ is finite) and $K\subseteq I$
is arbitrary, then $(X,Y)\mapstoo\proj_XY$ is continuous on the set
$\cD_{W_JW_K}^{J,K}$.
\end{Def}
The condition $\delta^*(X,Y)=W_JW_K$ means that $X$ and $Y$ are
\emph{almost opposite}, i.e. that there exist chambers $C\geq X$ and
$D\geq Y$ with $C\,\op\,D$.
\begin{Num}
\textbf{Remark}
For spherical buildings of rank 2 and for projective
spaces, this \emph{is} the common
notion of a topological building, cp.~Kramer \cite{KraDiss} and
K\"uhne \& L\"owen \cite{KueLoe}, K\"uhne \cite{Kue}.
\end{Num}
In the compact spherical case, there is a nice criterion
due to Grundh\"ofer and Van Maldeghem:
\begin{Prop}
\label{TheoProp}
If $\Delta$ is spherical, and if the topology on the vertex set
is compact, then $\Delta$ is a topological building if and only
if the chamber set is compact.

\proof
This follows from the closed graph theorem for maps
into compact spaces, see Grundh\"ofer \& Van Maldeghem \cite{GVM}.
\qed
\end{Prop}
We just mention the following result.
\begin{Prop}
\label{ResiduesProp}
Let $(\Delta^+,\Delta^-,\delta^*)$ be a topological twin
building, and let $X\in\Delta^+$ be spherical. Then $\Res(X)$ is
in a natural way a topological building. If $Z$ is another
simplex of the same type as $X$ 
(in either half of the twin building), then
$\Res(Z)$ is continuously isomorphic to $\Res(X)$.

\proof Pick $Y\in\Delta^-$ opposite $X$. Then $\proj$
induces an isomorphism $\Res(X)\cong\Res(Y)$. Let
$U,V\in\Res(X)$ be simplices with maximal $W_J$-distance.
We need to show that $\proj_VU$ is continuous. But
$\bar U=\proj_YU$ depends continuously on $U$, and
$\proj_VU=\proj_V\bar U$. For the second claim one uses the
following fact which is not difficult to prove (see the
proof by Tits \cite{Tits} 3.30 in the spherical case):
given simplices $X,X'$ (in the same half of the building),
there exists a simplex $X''$ in the other half which
is opposite both to $X$ and to $X'$,
\qed
\end{Prop}

\begin{Num}
\textbf{Example}
The building $\Delta(\CC^n)$, with the natural topology on its
vertex set $\Gr_1(\CC^n)\cup\Gr_2(\CC^n)\cup\cdots\cup\Gr_{n-1}(\CC^n)$
is a topological building by Proposition \ref{TheoProp},
since $\Fl(\CC^n)$ is compact.
\end{Num}
There is a natural topology on the group $\SL_n(\AA)$; the
evaluation map $\SL_n(\AA)\too \V_i$ induces a topology on the
set $\V_0\cup\cdots\cup\V_{n-1}$ of vertices of $\Delta^+(\AA^n)$ and
similarly on the vertices of $\Delta^-(\AA^n)$.
\begin{Thm}
\label{IsaTTB}
The twin building $(\Delta^+(\AA^n),\Delta^-(\AA^n),\delta^*)$
is a topological twin building.

\skop
Let $(C,D)\in\Delta^+(\AA^n)\times\Delta^-(\AA^n)$ be a pair of
opposite chambers, let $X\leq C$ be spherical of type $I\setminus J$
and $Y\leq D$ of type $I\setminus K$. Since $\SL_n(\AA)$
acts strongly transitively on the twin building, it acts
transitively on the set $\cD_{W_JW_K}^{J,K}$, and
$\cD_{W_JW_k}^{J,K}\cong\SL_n(\AA)/\SL_n(\AA)_{X,Y}$.
Now $\SL_n(\AA)_{X,Y}$ fixes $Z=\proj_XY$, and we have a continuous
map $\SL_n(\AA)/\SL_n(\AA)_{X,Y}\too\SL_n(\AA)/\SL_n(\AA)_Z$.
It follows that the map $(X',Y')\mapstoo \proj_{X'}Y'$ is continuous
on $\cD_{W_JW_K}^{J,K}$.
\qed
\end{Thm}
Let $\leqBruh$ denote the \emph{Bruhat order}
on $W/W_J$, for all $J\subseteq I$,
and put
\[
\cC_{\leqBruh wW_J}
(C_0)=\textstyle\bigcup\{\cC_{vW_J}(C_0)|\ vW_J\leqBruh wW_J\};
\]
this set is called the \emph{Schubert variety} corresponding to
$wW_J$. Let $m(wW_J)$ denote the $\ell$-length of the shortest coset
representative of $wW_J$.
\begin{Prop}
The Schubert varieties $\cC_{wW_J}(C_0)$ in $\Delta^+(\AA^n)$ are
CW-complexes, with Poincar\'e series
\[
\textstyle\sum_{vW_J\leqBruh wW_J}t^{2m(vW_J)}.
\]

\skop
Let $w$ be a shortest coset representative for $wW_J$, and let
$\cG_{s_1,\ldots,s_r}(C_0)$ denote the collection of all
(possibly stammering) galleries of type $(s_1,\ldots,s_r)$ (for
some fixed reduced expression $s_1\cdots s_r$ for $w$),
starting with the chamber $C_0$. It is
not difficult that to show that $\cG_{s_1,\ldots,s_r}(C_0)$ is
an iterated $\CP^1$-bundle
(sometimes called a \emph{Bott-Samelson cycle});
the total space is a smooth manifold of dimension $2m(wW_J)=2\ell(w)$.
These gallery spaces are also known as
\emph{Bott-Samelson desingularizations}
of Schubert varieties.
Consider the endpoint map
\[
\cG_{s_1,\ldots,s_r}(C_0)\too^\rho\cC_{\leqBruh wW_J}(C_0).
\]
The non-stammering galleries in $\cG_{s_1,\ldots,s_r}(C_0)$ are
mapped bijectively onto the Schubert cell\linebreak $\cC_{wW_J}(C_0)$. 
There is a canonical injection
\[
\cG_{s_1,\ldots,s_{r-1}}(C_0)\rInto\cG_{s_1,\ldots,s_r}(C_0)
\]
(by stammering at the end). The stammering galleries in
$\cG_{s_1,\ldots,s_r}(C_0)$ are either of this type, or galleries
which don't stammer at the end, but somewhere before the end.
Using this and the fact that
\[
\cG_{s_1,\ldots,s_r}(C_0)\too\cG_{s_1,\ldots,s_{r-1}}(C_0)
\]
is a $\CP^1$-bundle, one can use induction on $\ell(w)$ to
define a cellular map
\[
e^2\times e^{2\ell(w)-2}\too^\phi\cG_{s_1,\ldots,s_r}(C_0)
\]
which maps the boundary of the cell $e^2\times e^{2\ell(w)-2}$
onto the stammering galleries, see Mitchell \cite{Mitchell} Thm.~2.22 and
Kramer \cite{KraDiss} Sec.~4.1.
Now $\rho\circ\phi$ is an attaching map for a $2\ell(w)$-cell.
\qed
\end{Prop}
It follows that the based loop group $\Omalg\SU(n)$ has a
CW decomposition with Poincar\'e series
\[
\textstyle\sum_{wA_{n-1}\in\widetilde A_{n-1}/A_{n-1}}t^{2m(wA_{n-1})}
\]
where $W$ is the affine Weyl group of type $\widetilde A_{n-1}$
generated by $s_1,\ldots,s_n$,
and $A_{n-1}$ the subgroup generated by $s_1,\ldots,s_{n-1}$.

\subsection*{Knarr's construction for $\Delta^+(\AA^n)$}

We apply Knarr's construction to the the halves of the topological
twin building $(\Delta^+(\AA^n),\Delta^-(\AA^n),\delta^*)$.
The main ideas can be found in Mitchell's paper \cite{Mitchell},
although our approach (which is the same as Knarr's approach
\cite{Knarr}) is more based on geometric properties 
(i.e. the coordinatization of twin buildings), whereas Mitchell
makes strong use of the $BN$-pair. Let
$|\Delta^+(\AA^n)|$ denote the geometric realization of $\Delta^+(\AA^n)$.
By Theorem \ref{IsaTTB},
there is a canonical topology on the flag space $\Fl(\Delta^+(\AA^n))$
and we endow $|\Delta^+(\AA^n)|$ with the quotient topology induced
by the map
$\Fl(\Delta^+(\AA^n))\times|\simplex^{n-1}|\too|\Delta^+(\AA^n)|$.
The resulting space is denoted $|\Delta^+(\AA^n)|_\Knarr$.

More generally, assume that
$(\Delta^+,\Delta^-,\delta^*)$ is a topological twin building,
and that the panels are topological spheres. For example, in
$\Delta^+(\AA^n)$ the panels are homeomorphic to $\CP^1\cong\SS^2$.
Moreover, panels of the same type are homeomorphic by
Proposition \ref{ResiduesProp}. Let $m(s_i)$ denote the topological
dimension of the $i$-panels in $\Delta^+$.
It is proved in Kramer \cite{KraDiss} Prop.~2.0.2 that $m(s_i)=m(s_j)$
holds whenever $m_{ij}$ is odd. Thus we obtain a $\ZZ$-length
$m:W\too\ZZ$. For a Schubert cell $\cC_w(C_0)$, we have
\[
\cC_w(C_0)\cong\RR^{m(w)}.
\]
\begin{Prop}
(cp.~Knarr \cite{Knarr} and Mitchell \cite{Mitchell} 2.16)
Let $(\Delta^+,\Delta^-,\delta^*)$ be a topological twin building.
Assume that the panels are topological spheres. For $w\in W$ let
$X_w$ denote the image of $\cC_w(C_0)\times|\simplex^{r-1}|$
in $|\Delta^+|_\Knarr$.
Then for each $w\in W$, the set
$\bigcup\{X_v|\ v\prec w\}$ is contractible.

\skop
The proof is by induction on the length $\ell(w)$. Note that
$X_1=\{C_0\}\times|\simplex^{r-1}|$ is contractible.
Moreover, 
\[
\textstyle
\cC_{\preceq w}(C_0)/\bigcup\{\cC_u(C_0)|\ u\prec w\}\cong\SS^{m(w)},
\]
because this quotient is the one-point compactification of the
Schubert cell $\cC_w(C_0)$.

The inductive step is accomplished as follows.
First of all, $\bigcup\{X_v|\ v\prec u\}$ is contractible for all
$u\prec w$. Then it is not hard to see that 
$\bigcup\{X_v|\ v\preceq u\}/\bigcup\{X_v|\ v\prec u\}$ is also
contractible. This implies that $\bigcup\{X_v|\ v\preceq u\}$
is contractible. Finally, $\bigcup\{X_v|\ v\prec w\}$ is
homotopy equivalent to a wedge of such contractible spaces and
hence itself contractible.
\qed
\end{Prop}
\begin{Cor}
If $\Delta^+$ is spherical, then $|\Delta^+|_\Knarr$ is homeomorphic to
a sphere of dimension $m(w_0)+r-1$, where $w_0$ is the unique longest
element in the Coxeter group $W$, and $r$ is the rank of the building. 
In the non-spherical case, $|\Delta^+|_\Knarr$ is contractible.

\skop
In the spherical case, let $w_0\in W$ denote the longest element;
this is at the same time the unique maximal element in the Bruhat
order. Then $X_{w_0}\setminus\bigcup\{X_v|\ v\prec w_0\}
\cong\RR^{m(w_0)}\times\RR^{r-1}$, so
$X_{w_0}/\bigcup\{X_v|\ v\prec w_0\}\cong\SS^{m(w_0)}\wedge\SS^{r-1}$,
and this space is homotopy equivalent to $|\Delta^+|_\Knarr$,
since $\bigcup\{X_v|\ v\prec w_0\}$ is contractible.
Finally, it is not difficult to see that $X_{w_0}$ is a manifold, 
cp.~Kramer \cite{KraDiss} Prop.~4.2.1, hence
$|\Delta^+|_\Knarr$ is a compact manifold (of dimension at least
5) homotopy equivalent to a sphere, and thus to 
homeomorphic to $\SS^{m(w_0)+r-1}$
by the proof of the generalized Poincar\'e conjecture due to Smale
\cite{Smale}, Stallings \cite{Stallings} and Zeemann \cite{Zeemann}.
In the non-spherical case, $|\Delta^+|_\Knarr$ is a limit
of contractible spaces and hence itself contractible.
\qed
\end{Cor}
The Theorem above can be proved in much greater generality; it
suffices to assume that the panels are compact, connected, and
of finite covering dimension. Under these assumptions, the Corollary
holds up to homotopy equivalence, see Kramer \cite{KraDiss} Sec.~3.3.
On the other hand, if the vertices of the building are endowed
with the discrete topology, then the Theorem leads to a quick
proof of the Solomon-Tits Theorem.

\subsection*{Bott periodicity}
The crucial step in the proof of Bott periodicity is the
following observation. We have 
\[
\Res([E_0])\cong\Delta(\CC^n),
\]
and thus a natural map $\Gr_k(\CC^n)\rInto\V_k$. But the first
terms in the Poincar\'e series for these two spaces agree; the first
few shortest coset representatives for
$W/W_J$ and $W_K/W_{J\cap K}$ are the same, where
$J=\{2,\ldots,n\}$ and $K=\{2,\ldots,k-1,k+1,\ldots,n\}$.
For example, we have the following cell structure for $n=4$.
Put $I=\{1,2,3,4\}$ and let $J=\{2,3,4\}$ and $K=\{1,2,4\}$;
thus we have the following Coxeter groups.
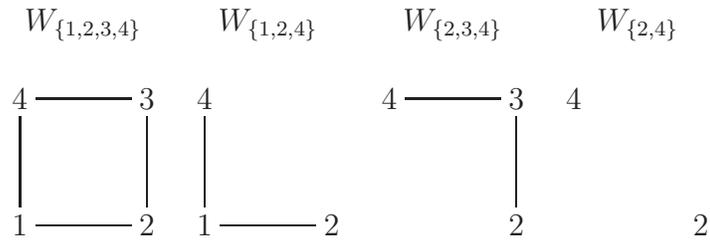
\[
\begin{array}{cccc}
W_{\{1,2,3,4\}} & W_{\{1,2,4\}} & W_{\{2,3,4\}} & W_{\{2,4\}} \\ \\
{\xymatrix@!=2em{4 \ar@{-}[r]\ar@{-}[d] & 3 \ar@{-}[d]\\
1\ar@{-}[r] & 2}} &
{\xymatrix@!=2em{4 \ar@{-}[d] \\ 1\ar@{-}[r] & 2}} &
{\xymatrix@!=2em{4 \ar@{-}[r] & 3\ar@{-}[d] \\ & 2}} &
{\xymatrix@!=2em{4  \\ & 2}}
\end{array}
\]
For $W_{\{1,2,4\}}/W_{\{4,2\}}$ we have the following
Bruhat order for the shortest coset representatives,
cp.~Scharlau \cite{Scharlau} 2.5.
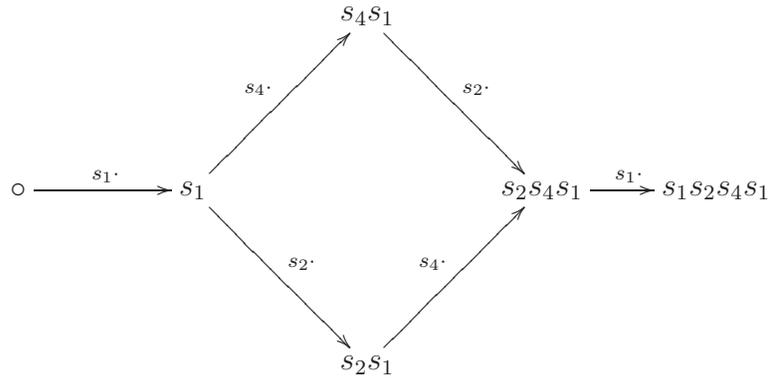
\[
\xymatrix@!=3.5em{
&& s_4s_1 \ar[dr]^{s_2\cdot} \\
\circ \ar[r]^{s_1\cdot} & s_1 \ar[ur]^{s_4\cdot}\ar[dr]^{s_2\cdot}
&& s_2s_4s_1\ar[r]^{s_1\cdot} & s_1s_2s_4s_1 \\
&& s_2s_1 \ar[ur]^{s_4\cdot}
}
\]
Thus we see that $\Gr_2(\CC^4)$ has a cell decomposition as
$e^0\cup e^2\cup e^4\cup e^4\cup e^6\cup e^8$.
The Bruhat order for the shortest coset representatives of
$W_{\{1,2,3,4\}}/W_{\{2,3,4\}}$ starts as
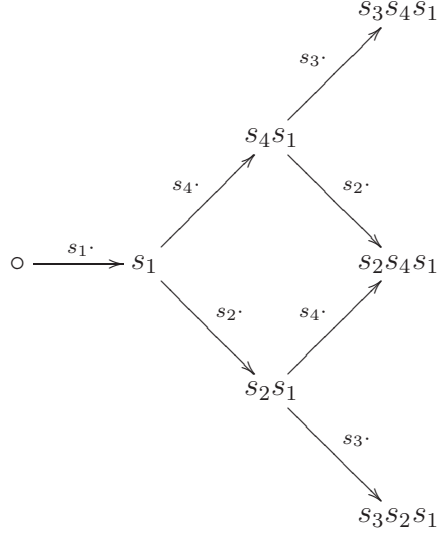
\[
\xymatrix@!=2em{
&&& s_3s_4s_1 \\
&& s_4s_1 \ar[ur]^{s_3\cdot}\ar[dr]^{s_2\cdot} \\
\circ \ar[r]^{s_1\cdot} & s_1 \ar[ur]^{s_4\cdot}\ar[dr]^{s_2\cdot}
&& s_2s_4s_1 & {}\\
&& s_2s_1 \ar[ur]^{s_4\cdot}\ar[dr]^{s_3\cdot}\\
&&& s_3s_2s_1 
}
\]
(and continues infinitely to the right).
Accordingly, the cell structure of $\V_2\cong\Omalg\SU(4)$ is
$e^0\cup e^2\cup e^4\cup e^4\cup e^6\cup e^6\cup e^6\cup\cdots$;
in particular, the $5$-skeleton of $\V_2$ is the same as for
$\Gr_2(\CC^4)$. The general result is as follows.
\begin{Prop}
Suppose that $n=2k$ is even. Then the inclusion
$\Gr_k(\CC^{2k})\rInto\V_k$ is a $2k$-equivalence.
\qed
\end{Prop}
The proof for Bott periodicity is now as follows. If $k$ is large,
then $\Gr_k(\CC^{2k})$ is a good approximation for the
classifying space $\mathrm{BU}$; so in the limit, we obtain
a homotopy equivalence
$\mathrm{BU}\simeq\Omega_\alg\SU$. There is one problem, though;
we have only considered the space $\Omalg\SU(n)$ consisting of all
based loops which can be expressed as \emph{Laurent polynomials},
whereas topologists consider the space $\Omega_{\mathrm{cts}}\SU(n)$ of
all \emph{continuous} (or smooth)
based loops. So the proof is not yet finished.
At this point it is convenient to introduce Quillen's 
\emph{space of special paths}
\[
\cS=\left\{\left.\left(
t\mapstoo g(e^{2\pi\ti t})e^{2\pi\ti tX}g^\#(1)\right)\right| \ 
X\in\fc,\,g\in L_\alg\SU(n)\right\}\subseteq (\SU(n),\one)^{([0,1],0)},
\]
where $\fc\subseteq\fp_n$ is a Weyl chamber. The group $L_\alg\SU(n)$
acts in a natural way on $\cS$. But $\fc\subseteq\fp_n\subseteq\fX$
can also be identified with the image of the
chamber $\{[E_0],\ldots,[E_{n-1}]\}$ under the Veronese representation.
This identification extends in a natural way to an
$L_\alg\SU(n)$-equivariant homeomorphism 
\[
\cS\too^\cong|\Delta^+(\AA^n)|_\Knarr\subseteq\fX
\]
as follows: a path $\gamma\in\cS$ is mapped to its
'logarithmic derivative'
$\frac{1}{2\pi\ti}\left(\frac{d}{dt}\gamma(t)\right)\gamma(t)^{-1}$
which is a path in the tangent space $T_1\SU(n)$. (Conversely,
a smooth path in the tangent space can be read as a differential
equation whose solution starting at $1$ is a path in the Lie group.)
If we put $z=e^{2\pi\ti t}$, then
\begin{align*}
\textstyle
\frac{d}{dt}\left(g(z)e^{2\pi\ti tX}g^\#(1)\right) &= \textstyle
2\pi\ti
\left((\zpz g(z))e^{2\pi\ti tX}g^\#(1)+g(z)Xe^{2\pi\ti tX}g^\#(1)\right) \\
\textstyle
\left(\frac{d}{dt}\left(g(z)e^{2\pi\ti tX}g^\#(1)\right)\right)
\left(g(z)e^{2\pi\ti tX}g^\#(1) \right)^{-1} &= \textstyle
2\pi\ti \left(\left(\zpz g(z)\right)g^\#(z)+g(z)Xg^\#(z)\right).
\end{align*}
\begin{Lem}
The map 
\[
\textstyle\gamma\mapstoo
\frac{1}{2\pi\ti}\left(\frac{d}{dt}\gamma(t)\right)\gamma(t)^{-1}
\]
is an $L_\alg\SU(n)$-equivariant homeomorphism
$\cS\too|\Delta^+(\AA^n)|_\Knarr\subseteq\fX$, where the action
on $\Delta^+(\AA^n)$ is the standard one, and the action on
Quillen's space of special paths is by
$\gamma(t)\mapstoo g(e^{2\pi\ti t})\gamma(t) g^{-1}(1)$.
\qed
\end{Lem}
In particular, $\cS$ is contractible. The endpoint map
$(\SU(n),\one)^{([0,1],0)}\too\SU(n)$ yields the universal bundle
\[
\Omega_{\mathrm{cts}}\SU(n)\too
(\SU(n),\one)^{([0,1],0)}\too\SU(n),
\]
(where $\mathrm{cts}$ refers to the group of continuous loops)
and as a subbundle we have
\[
\Omalg\SU(n)\too\cS\too\SU(n).
\]
Since the total spaces of both bundles are contractible,
the inclusion $\Omalg\SU(n)\rInto\Omega_{\mathrm{cts}}\SU(n)$ is a 
weak (and therefore also a strong) homotopy equivalence.
We have proved Quillen's following result.

\begin{Thm}(Quillen)
The orbit space $|\Delta^+(\AA^n)|_\Knarr/\Omalg\SU(n)$ is
homeomorphic to $\SU(n)$.
\qed
\end{Thm}
Recall that $\Gr_k(\CC^n)$ is a good approximation of the classifying
space $\mathrm{BU}$ in small dimensions.
\begin{Cor}\textbf{\em(Unitary Bott Periodicity)}
The inclusion
\[
\Gr_k(\CC^{2k})\too\Omega_\alg\SU(k)\rInto^\simeq
\Omega_{\mathrm{cts}}\SU(k)
\]
is a $2k$-equivalence. In the limit, the natural map
\[
\mathrm{BU}\too\Omega_{\mathrm{cts}}\SU
\]
is a homotopy equivalence.
\qed
\end{Cor}
This implies in particular the famous Bott isomorphisms
$\pi_{2k}(\U)=0$ and $\pi_{2k+1}(\U)\cong\ZZ$, for all $k\geq 0$.
Note that we have also proved that the $\Omalg\SU(n)$-orbit space map
\[
|\Delta^+(\AA^n)|_\Knarr\too\SU(n)
\]
is a universal classifying bundle for $\Omalg\SU(n)$.

\section{Real forms and compact symmetric spaces}
So far, we have discussed the group $\SL_n(\AA)$ which is the proper
generalization of the complex group $\SL_n(\CC)$. In this last section
we consider briefly how real groups fit into the picture;
more details can be found in Mitchell \cite{Mitchell} and
in Gro\ss, Heintze, Kramer \& M\"uhlherr \cite{GHKM}.
Consider the involution $\iota$ given by
\[
\sum_\fin X_\nu z^\nu\mapstoo \sum_\fin \bar X_\nu z^\nu.
\]
The group of $\iota$-fixed elements in $\SL_n(\AA)$ is
$\SL_n(\RR[z,1/z])$, and there is a corresponding twin building
over $\RR[z,1/z]$ which is defined exactly in the same way as the one
over $\CC[z,1/z]$ considered so far.
Note however that we cannot interpret the elements of
$\SL_n(\RR[z,1/z])$ as loops in $\SL_n(\RR)$. Instead, we view the
elements of $\SL_n(\RR[z,1/z])$ as paths
\begin{align*}
[0,1]&\too L_\alg\SL_n(\CC)\\
t &\mapstoo g(e^{\ti\pi t}).
\end{align*}
These paths
have the special property that they start and end in $\SL_n(\RR)$.
If we intersect $\SL_n(\RR[z,z^{-1}])$ with $L_\alg\SU(n)$, then we
obtain the group
\[
\SL_n(\RR[z,1/z])\cap L_\alg\SU(n)
\]
consisting of paths in $\SU(n)$ which start and end in $\SO(n)$.
Similarly as before, this group is homotopy equivalent with the based
loop space
\[
\Omega_{\mathrm{cts}}(\SU(n)/\SO(n)).
\]
These loop spaces of compact Riemannian symmetric spaces play an
important r\^ole in topology.
They can be used to prove the other versions of Bott periodicity
(real and quaternionic), see Mitchell \cite{Mitchell}.

\smallskip
\textbf{Acknowledgement}
I am indebted to Peter Abramenko for sharing some of his insights,
and to Theo Grundh\"ofer for some remarks on the paper.
I would particularly like to thank Bernhard M\"uhlherr, who
introduced me to twin buildings in Oberwolfach eight years ago
and convinced me of their usefulness and beauty;
many ideas in this paper stem from discussions with him during the
last years.



\begin{thebibliography}{mm}

\bibitem{AbraLNM}
P. Abramenko,
\emph{Twin buildings and applications to $S$-arithmetic groups},
Springer LNM 1641 (1996).

\bibitem{AM}
P. Abramenko and B. M\"uhlherr,
Pr\'esentations de certaines $BN$-paires jumel\'ees
comme sommes amalgam\'ees,
C. R. Acad. Sci. Paris 325 (1997) 701--706.

\bibitem{AR}
P. Abramenko and M. Ronan,
A characterization of twin buildings by twin apartments,
Geom. Dedicata 73 (1998) 1--9.

\bibitem{AVM}
P. Abramenko and H. Van Maldeghem,
1-twinnings of buildings,
to appear in Math. Z. (2001).

\bibitem{Rou}
V. Back-Valente, N. Bardy-Panse,
H. Ben Massaoud, and G. Rousseau,
Formes presque-d\'eploy\'ees des
alg\`ebres de Kac-Moody: Classification et racines
relatives,
J. Algebra 171 (1995) 43--96.

\bibitem{MueCox2}
N. Brady, J.P. McCammond, B. M\"uhlherr, and W.D. Neumann,
Non-rigid Coxeter and Artin groups,
Preprint, Dortmund (2000).

\bibitem{BridHaef}
M.R. Bridson and A. Haefliger,
\emph{Metric spaces of non-positive curvature},
Grundlehren der Mathema\-ti\-schen Wissenschaften 319,
Springer-Verlag, Berlin (1999).

\bibitem{Brown}
K. Brown,
{\em Buildings},
Springer Verlag (1989).

\bibitem{BurnsSpatzier}
Burns, K., Spatzier, R.:
{On topological Tits buildings and their classification},
Publ. Math. I.H.E.S. 65 (1987) 5--34.

\bibitem{CD}
R. Charney and M. Davis,
When is a Coxeter system determined by its Coxeter group?,
J. London Math. Soc. 61 (2000) 441--461.

\bibitem{DS}
A. Dress and R. Scharlau,
Gated sets in metric spaces,
Aequationes math. 34 (1987) 112--120.

\bibitem{GHKM}
C. Gro\ss, E. Heintze, L. Kramer, and B. M\"uhlherr,
manuscript in preparation.

\bibitem{GK}
T. Grundh\"ofer and N. Knarr,
Topology in generalized quadrangles,
Topology Appl. 34 (1990) 139--152.

\bibitem{Flaghom1}
T. Grundh\"ofer, N. Knarr, and L. Kramer,
Flag-homogeneous compact connected polygons,
Geom. Dedicata 55 (1995) 95--114.

\bibitem{Flaghom2}
T. Grundh\"ofer, N. Knarr, and L. Kramer,
Flag-homogeneous compact connected polygons II,
Geom. Dedicata 83 (2000) 1--29.

\bibitem{Flaghom3}
T. Grundh\"ofer, N. Knarr, and L. Kramer,
The classification of compact homogeneous buildings,
manuscript in preparation.

\bibitem{GVM}
T. Grundh\"ofer and H. Van Maldeghem,
Topological polygons and affine buildings of rank three,
Atti Sem. Mat. Fis. Univ. Modena 38 (1990) 459--479.

\bibitem{HeintzeLiu}
E. Heintze and X. Liu,
Homogeneity of infinite-dimensional isoparametric submanifolds,
Ann. Math. 149 (1999) 149--181.

\bibitem{HPTT}
E. Heintze, R. Palais, C.-l. Terng, and G. Thorbergsson,
Hyperpolar actions on symmetric spaces, in:
\emph{Geometry, Topology, \& Physics, for Raoul Bott},
Cambridge, MA, 1993. 
S.-T. Yau ed.,
International Press (1995) 214--245.

\bibitem{Jae}
M. J\"ager,
\emph{Topologische Geb\"aude},
Dissertation, Univ. Kiel (1994)

\bibitem{Knarr}
N. Knarr,
The nonexistence of certain topological polygons,
Forum Math. 2 (1990) 603--612.

\bibitem{KK}
N. Knarr and L. Kramer,
Projective planes and isoparametric hypersurfaces,
Geom. Dedicata 58 (1995) 193--202. 

\bibitem{KraDiss}
L. Kramer,
{\em Compact polygons},
Dissertation, Univ. T\"ubingen (1994).\\
available as math.DG/0104064 in the Mathematics ArXiv,\\
\texttt{\small http://front.math.ucdavis.edu/math.DG/0104064}

\bibitem{KraHabil}
L. Kramer,
\emph{Homogeneous spaces, Tits buildings, and isoparametric hypersurfaces},
to appear in Mem. Amer. Math. Soc. (2001).

\bibitem{Kue}
R. K\"uhne,
\emph{Topologische sph\"arische Tits-Geb\"aude},
Dissertation, Univ. Braunschweig (1994)\\
Mitt. Math. Seminar Giessen 219 (1994).

\bibitem{KueLoe}
R. K\"uhne and R. L\"owen,
Topological projective spaces,
Abh. Math. Seminar Univ. Hamburg 62 (1992) 1--9.

\bibitem{Mitchell}
S.A. Mitchell,
Quillen's theorem on buildings and the loops on a 
symmetric space,
L'enseignement math\'ematique 34 (1988) 123--166.

\bibitem{Mu3}
B. M\"uhlherr,
A rank 2 characterization of twinnings,
European J. Comb. 10 (1998) 603--612.

\bibitem{Mu2}
B. M\"uhlherr,
\emph{On the existence of 2-spherical twin buildings},
Habilitationsschrift, Univ. Dortmund (1999).

\bibitem{MueCox1}
B. M\"uhlherr,
On isomorphisms between Coxeter groups,
Designs, Codes and Cryptography 21 (2000) 188--189.

\bibitem{MR}
B. M\"uhlherr and M. Ronan,
Local to global structure in twin buildings,
Invent. Math. 122 (1995) 71--81.

\bibitem{PalaisTerng}
R. Palais and C.-l. Terng,
\emph{Critical point theory and submanifold geometry},
Springer LNM 1353 (1988).

\bibitem{PT}
U. Pinkall and G. Thorbergsson,
Examples of infinite dimensional
isoparametric submanifolds,
Math. Z. 205 (1990) 279--286.

\bibitem{PresSeg}
A. Pressley and G. Segal,
{\em Loop groups},
Oxford University Press (1986), corr.~reprint (1988).

\bibitem{Ronan}
M. Ronan,
\emph{Lectures on buildings},
Perspectives in Mathematics 7,
Academic Press, Boston, MA (1989).

\bibitem{Ron2}
M. Ronan,
Local isometries of twin buildings,
Math. Z. 234 (2000) 435--455.

\bibitem{RT}
M. Ronan and J. Tits,
Twin trees I,
Invent. Math. 116 (1994) 463--479.

\bibitem{Rou2}
G. Rousseau, 
On forms of Kac-Moody algebras, in:
\emph{Algebraic groups and their generalizations: quantum
and infinite-dimensional methods},
University Park, PA, 1991,
W. Haboush and B. Parshall ed.,
Proc. Sympos. Pure Math. 56 (1994) 393--399.

\bibitem{Salz}
H. Salzmann,
Topological planes,
Adv. Math. 2 (1967) 1--60.

\bibitem{CPP}
H. Salzmann, D. Betten, T. Grundh\"ofer, H. H\"ahl, R. L\"owen, and
M. Stroppel,
{\em Compact projective planes},
De Gruyter Expos. in Math. 21, Berlin (1995).

\bibitem{Scharlau}
R. Scharlau,
Buildings,
in: \emph{Handbook of incidence geometry}, 
F. Buekenhout ed.,
North Holland, Amsterdam (1995) 477--645.

\bibitem{Schroth}
A. Schroth,
\emph{Topological circle planes and topological quadrangles},
Pitman RNM 337 (1995).

\bibitem{Smale}
S. Smale,
Generalized Poincar\'e's conjecture in dimensions greater
than four,
Ann. Math. 74 (1961) 391--406.

\bibitem{Stallings}
J. Stallings,
Polyhedral homotopy-spheres,
Bull Amer. Math. Soc. 66 (1960) 485--488.

\bibitem{Te}
C.-l. Terng,
Recent progress in submanifold geometry, in:
\emph{AMS Summer Research Institute on Differential Geometry},
Los Angeles, CA, 1990,
R. Green and S.-T. Yau ed.,
Proc. Symp. Pure Math. 54 (1993) 439--484.

\bibitem{Thorbg}
G. Thorbergsson,
Isoparametric foliations and their buildings,
Ann. Math. 133 (1991) 429--446.

\bibitem{ThorbgHandbook}
G. Thorbergsson,
A survey on isoparametric hypersurfaces and their generalizations,
in: \emph{Handbook of differential geometry},
F. Dillen and L. Verstraelen ed.,
North Holland, Amsterdam (2000) 963--995.

\bibitem{Tits}
J. Tits,
\emph{Buildings of spherical type and finite $BN$-pairs},
Springer LNM 386 (1974).

\bibitem{StAndrews}
J. Tits,
Buildings and group amalgamations,
in: \emph{Proceedings of groups -- St Andrews 1985},
Ed. E.F. Robertson and C.M. Campbell,
Cambridge University Press, Cambridge (1986) 110--127.

\bibitem{TitsRes}
J. Tits,
\emph{Immeubles jumel\'es},
R\'esum\'e de cours.\\
Annuaire du Coll\`ege de France, 89$^e$ ann\'ee (1988--1989) 81--95.\\
Annuaire du Coll\`ege de France, 90$^e$ ann\'ee (1989--1990) 87--103.\\
Annuaire du Coll\`ege de France, 98$^e$ ann\'ee (1997--1998) 97--112.

\bibitem{Durham}
J. Tits,
Twin buildings and groups of Kac-Moody type,
in: \emph{Groups, combinatorics \& geometry},
Durham, 1990,
M. Liebeck and J. Saxl ed.,
London Math. Soc. Lecture Notes,
Cambridge (1992) 249--286.

\bibitem{HVM}
H. Van Maldeghem,
\emph{Generalized polygons},
Birkh\"auser Monographs in Math. 93, Basel (1998).

\bibitem{Zeemann}
E.C. Zeeman,
The generalized Poincar\'e conjecture,
Bull Amer. Math. Soc. 67 (1961) 270.

\end{thebibliography}
\end{document}